\documentclass[a4paper, 11pt, leqno]{amsart}
\usepackage[utf8]{inputenc}
\usepackage[usenames,dvipsnames]{color}
\usepackage{amssymb, amsmath, latexsym, graphics, graphicx,ytableau}
\usepackage[hidelinks]{hyperref}
\usepackage{amsfonts}
\usepackage{amsthm}
\usepackage{enumerate}

\makeatletter
\newcommand{\leqnomode}{\tagsleft@true}
\newcommand{\reqnomode}{\tagsleft@false}
\makeatother

\newtheorem{theorem}{Theorem}[section]
\newtheorem{corollary}[theorem]{Corollary}
\newtheorem{lemma}[theorem]{Lemma}
\newtheorem{proposition}[theorem]{Proposition}
\newtheorem{question}[theorem]{Question}

\newtheorem{remark}[theorem]{Remark}
\theoremstyle{definition}

\newcommand{\B}{\mathbb{B}}
\newcommand{\N}{\mathbb{N}}
\newcommand{\R}{\mathbb{R}}

\newcommand{\T}{\mathbb{R}_\mathrm{max}}
\newcommand{\Q}{\mathbb{Q}}
\newcommand{\Z}{\mathbb{Z}}

\newcommand{\cR}{\mathcal{R}}

\newcommand{\cJ}{\mathrel{\mathcal{J}}}
\newcommand{\cI}{\mathcal{I}}
\newcommand{\cG}{\mathcal{G}}
\newcommand{\cM}{\mathcal{M}}
\newcommand{\cS}{\mathcal{S}}

\newcommand{\cX}{\mathcal{X}}
\newcommand{\cZ}{\mathcal{Z}}
\newcommand{\Qmax}{\Q_\mathrm{max}}
\newcommand{\Zmax}{\Z_\mathrm{max}}

\newcommand{\floor}[1]{\left \lfloor #1 \right \rfloor }

\newcommand{\lanran}[1]{\left\langle #1 \right\rangle}

\title[Generating Sets, Presentations, and Growth]{Generating Sets, Presentations, and Growth of Tropical Matrix Monoids}
\begin{document}
\date{\today}
\maketitle
\begin{center}
THOMAS AIRD\footnote{Email \texttt{Thomas.Aird@manchester.ac.uk}.} \\ \ \\
Department of Mathematics, University of Manchester, \\
Manchester M13 9PL, UK. \\
Heilbronn Institute for Mathematical Research, Bristol, UK.
\end{center}

\begin{abstract}
We construct minimal and irredundant generating sets for a family of submonoids of the monoid of $n \times n$ upper triangular matrices over a commutative semiring.
We show that the monoid of $n \times n$ matrices over the tropical integers, $M_n(\Zmax)$, is finitely generated if and only if $n \leq 2$, and finitely presented if and only if $n = 1$. 
Minimal and irredundant generating sets are explicitly constructed when $n \leq 3$. 
We then construct a presentation for the monoid of $n \times n$ upper triangular matrices over the tropical integers, $UT_n(\Zmax)$, demonstrating that it is finitely presented for all $n \in \N$. 
Finally, we establish upper bounds on the polynomial degree of the growth function of finitely generated subsemigroups of the monoid of $n \times n$ matrices over a bipotent semiring and show that these bounds are sharp for the tropical semiring.
\end{abstract}

\section{Introduction}
\par Constructing minimal and irredundant generating sets for semigroups is a widely studied area of research (see, for example, \cite{SingularTransformation,0SimpleSemigroups}) and is related to the classical problem of calculating the \emph{rank} of a semigroup, that is, the minimum cardinality of a generating set of a semigroup \cite{GpRank,SgpRank}. 

\par The monoid of $n \times n$ matrices over the tropical semiring, $M_n(\T)$, has attracted considerable interest. This is due to its many useful properties, including its ability to admit faithful representations of semigroups which cannot be faithfully represented by matrices over fields \cite{ARStylic,CJKMPlacticLikeRep,2by2,JKPlactic}. 

\par Recently, there has been research into constructing minimal generating sets for matrix monoids.
In particular, East, Jonu\v sas and Mitchell \cite{Gen2by2} found generating sets for $2 \times 2$ full matrix monoids over the min-plus natural number semiring, max-plus natural number semiring, and their finite quotients. 
These generating sets were later shown to be minimal by Hivert, Mitchell, Smith, and Wilson \cite{MinGens} who further found minimal generating sets for a number of submonoids of the monoid of $n \times n$ boolean matrices.

\par Beyond generating sets, another interesting problem in semigroup theory is determining whether a semigroup admits a finite presentation, that is, whether a semigroup has a finite generating set in which all equalities can be deduced from a finite set of relations. 
This is a very active area of research, with new presentations being constructed for many semigroups which naturally embed into $M_n(\T)$, \cite{EPresentationPK,EHOOPNumericalSemigroups,MWPresentations}.

\par Relating to both generating sets and presentations is growth. The growth rate of a semigroup is an important invariant in geometric semigroup theory. Understanding the growth of a semigroup provides information about the geometry and structure of the semigroup, \cite{Gray2017}. For instance, Gromov's theorem on groups of polynomial growth states that a finitely generated group has polynomial growth if and only if it has a nilpotent subgroup of finite index \cite{Gromov}.

\par The growth rate of subsemigroups of $M_n(\T)$ was first studied by d'Alessandro and Pasku in \cite{CombinatorialSemigroups}. In particular, they investigated the growth of finitely generated subsemigroups of $M_n(S)$, when $S$ is a commutative bipotent semiring, and show that for any finitely generated subsemigroup, the growth function is bounded above by a polynomial. However, the degree of the polynomial is dependent on the dimension of the matrices and the number of unique entries in the matrices in the generating set. As a result, different generating sets for the same semigroup can give different upper bounds on the polynomial degree of the growth rate.

\par In this paper, we study matrix monoids with a focus on these three important properties: generating sets, presentations, and growth rate. Moreover, we are interested in how these properties apply to the monoid of $n \times n$ matrices over the tropical integer semiring, $\Zmax$.

\par In particular, this paper comprises 6 sections, including this introduction. In Section \ref{sec:Prelim}, we introduce notation and definitions used throughout.
In Section \ref{UTGenSec}, we describe the minimal and irredundant generating sets of the monoid of upper triangular matrices over a commutative semiring with the diagonal entries coming from a fixed submonoid of the multiplicative semigroup of the semiring. 
Thus, showing that the monoid of $n \times n$ upper triangular matrices over the tropical integers, $UT_n(\Zmax)$, is finitely generated for all $n \in \N$ and the monoid of $n \times n$ unitriangular matrices over $\Zmax$ is finitely generated if and only if $n = 1$. 
\par In Section \ref{FullGenSec}, we turn our attention to full matrix monoids, constructing a finite minimal generating set for $M_2(\Zmax)$. For $n \geq 3$, we show that if $S$ is an anti-negative semifield, then $M_n(S)$ is finitely generated if and only if $S$ is finite. Thus, $M_n(\Zmax)$ is finitely generated if and only if $n \leq 2$. 
We then construct a minimal and irredundant generating set for $M_3(\Zmax)$ and show that the subsemigroup of $M_3(\Zmax)$ generated by all regular matrices can be generated by four elements.
\par In Section \ref{PresentationSec}, we show that $UT_n(\Zmax)$ is finitely presented, for all $n \in \N$, by showing that every word over the generators can be rewritten into a normal form. We then use this presentation to give a different finite presentation for $UT_n(\Zmax)$ using the minimal and irredundant generating set found in Section~\ref{UTGenSec}.
\par In Section~\ref{GrowthSec}, we find, for a commutative bipotent semiring $S$, an upper bound on the polynomial degree for the growth function of any finitely generated subsemigroup of $M_n(S)$, producing a similar bound for $UT_n(S)$.
Moreover, when $S = \T$, we get bounds dependent only on $n$ and the rank of the free abelian subgroup which the finite entries in the matrices of the generating set generate as a group, which is independent of the generating set. 
Finally, we show that these bounds are sharp for $M_n(\T)$ and $UT_n(\T)$ for all $n \in \N$, by giving examples of finitely generated subsemigroups of $M_n(\T)$ and $UT_n(\T)$ which attain these bounds.

\section{Preliminaries} \label{sec:Prelim}
\par Let $\N$ denote the set of positive integers and $\N_0$ be the set of non-negative integers. For a semigroup $\cS$, $X \subseteq \cS$ is a \emph{(semigroup) generating set} for $\cS$, if $\cS$ is the smallest subsemigroup of $\cS$ containing $X$, in this case, we write $\lanran{X} = \cS$.
For a group $\cG$, $X$ is a \emph{group generating set} for $\cG$ if $X \cup X^{-1} \cup \{1_\cG\}$ is a semigroup generating set for $\cG$. 
We say a generating set $X$ for $\cS$ is \emph{minimal} if $|X| \leq |Y|$ for any generating set $Y$ for $\cS$ and say an element $x \in X$ is \emph{irredundant} if $X \setminus \{x\}$ is not a generating set for $\cS$. If every $x \in X$ is irredundant then we say $X$ is \emph{irredundant}. More generally, we say a set $X$ is \emph{minimal with a given property} if it has the property and $|X| \leq |Y|$ for any set $Y$ that has the property, and say a set is \emph{irredundant with a given property} if it has the property and no proper subset of it has the property.

\par We call $x \in \cS$ a \emph{unit} of a monoid $\cS$ if there exists $x^{-1} \in \cS$ such that $xx^{-1} = x^{-1}x = 1_\cS$. Let $U(\cS)$ be the \emph{group of units} of $\cS$, that is, the set of all units in $S$. We say a non-unit $x \in \cS$ is \emph{prime} if, for every product $x = uv$, exactly one of $u$ or $v$ is a unit.
For a monoid $\cS$, we define Green's $\cJ$-relation to be the equivalence relation on $\cS$ defined by $x \cJ y$ if and only if $\cS x\cS = \cS y\cS$. For $a \in \cS$, denote the $\cJ$-class containing $a$ by $J_a$. We call $J$ a \emph{prime} $\cJ$-class if every element of $J$ is prime.
As multiplying by a unit keeps elements in the same $\cJ$-class, it is easy to see that every generating set of $\cS$ contains a representative from each prime $\cJ$-class of $\cS$. 

\par Let $S$ be a (unital) semiring, that is, a set $S$ with two binary operations $+$ and $\cdot$ such that multiplication distributes over addition, $(S,+)$ is a commutative monoid with identity $0_S$, and $(S,\cdot)$ is a monoid with identity $1_S$ such that
$x0_S = 0_Sx = 0_S$ for all $x \in S$. We say $S$ is \emph{commutative} if $(S, \cdot)$ is commutative, and a \emph{semifield} if $(S^*, \cdot)$ is an abelian group where $S^* = S \setminus \{0_S\}$.

\par For a semiring $S$, let $U(S)$ be the \emph{group of units} of $(S,\cdot)$. We say $x \in S$ is \emph{additively invertible} if there exists $y \in S$ such that $x + y = 0_S$. Let $V(S)$ be the subset of additively invertible elements of $S$, i.e. the group of units of $(S,+)$. 
Note that, if $x,y \in V(S)$ and $z \in S$, then $x+y \in V(S)$ and $zx,xz \in V(S)$. Thus, $V(S)$ is a (possibly non-unital) ring and $V(S) = S$ if and only if $1_S \in V(S)$. 
We say $S$ is \emph{anti-negative} if for $x,y \in S$, $x + y = 0_S$ if and only if $x = y = 0_S$, that is, if $V(S) = \{0_S\}$. 

\par Let $M_n(S)$ be the monoid of all $n \times n$ matrices with entries in $S$ under matrix multiplication and $UT_n(S)$ be the submonoid of all $n \times n$ upper triangular matrices over $S$, that is, matrices with $0_S$ entries below the diagonal. Then, for a fixed submonoid $T$ of $(S,\cdot)$, let $UT_n^T(S)$ be the submonoid of $UT_n(S)$ in which all diagonal entries are from $T$. Note that $UT_n^S(S) = UT_n(S)$ and $UT_n^{\{1_S\}}(S)$ is the monoid of all $n \times n$ unitriangular matrices over $S$.

\par Finally, we define some matrices we use throughout. For $1 \leq i \leq n$, we let $A_i(\lambda) \in UT_n(S)$ be the diagonal matrix with $1_S$ on the diagonal apart from $\lambda$ as the $(i,i)$th entry, and, for $1 \leq i < j \leq n$, let $E_{ij}(\lambda) \in UT_n(S)$ be the matrix where all diagonal entries are $1_S$, $(E_{ij})_{ij} = \lambda$, and all other entries are $0_S$. We sometimes write $E_{ij}$ to denote $E_{ij}(1_S)$.
\par Let $\B = \{0,1\}$ with addition and multiplication given by maximum and minimum respectively be the \emph{boolean semifield}, and let $\T = \R \cup \{-\infty\}$ with addition and multiplication given by maximum and addition respectively be the \emph{tropical semiring}.
We denote the subsemirings of tropical integers and tropical rationals by $\Zmax = \T \cap (\Z \cup \{-\infty\})$ and $\Qmax = \T \cap (\Q \cup \{-\infty\})$ respectively. 
\par We begin, by introducing two lemmas which we require for the following two sections.
\begin{lemma} \label{lem:CommJClass}
Let $S$ be a commutative semiring. Then, $xy$ is a unit if and only if $x$ and $y$ are units.
\end{lemma}
\begin{lemma} \label{Jclass}
Let $S$ be a commutative semiring and $X \in \cS$ where $\cS = M_n(S)$ or  $UT_n^T(S)$ for some $T$ a submonoid of $(S,\cdot)$. If $X \cJ I_n$ in $\cS$, then $X$ is a unit in $\cS$.
\end{lemma}
\begin{proof}
If $X \cJ I_n$, then there exists $A,B \in \cS$ such that $AXB = I_n$. Hence, by the main theorem in  \cite{MatrixJClass}, $XBA = BAX = I_n$, and thus $X \in U(S)$.
\end{proof}

\section{Generating sets for upper triangular matrix monoids} \label{UTGenSec}
In this section, we produce minimal and irredundant generating sets for $UT_n^T(S)$ when $S$ is a commutative semiring and $T$ is any submonoid of $(S,\cdot)$. 
By choosing different $T$ be obtain many interesting submonoid of $UT_n(S)$, in particular, there has been a lot of interest in the submonoids when $T = S$, $\{1_S\}$, $\{1_S,0_S\}$, or $U(S)$, \cite{MatrixOverSemirings, HZLEquationTheoriesUniZero, IJKTropicalGroups}. 

\par We begin by characterising exactly when a matrix in $UT_n^T(S)$ is invertible in $UT_n^T(S)$.

\begin{lemma} \label{InvCond}
Let $n \in \N$, $S$ be a commutative semiring, and $T$ be a submonoid of $(S,\cdot)$. Then, $X \in UT_n^T(S)$ is invertible in $UT_n^T(S)$ if and only if $X_{ii} \in U(T)$ for $1 \leq i \leq n$ and $X_{ij} \in V(S)$ for $1 \leq i < j \leq n$.
\end{lemma}
\begin{proof}
\par Let $X \in UT_n^T(S)$ be invertible in $UT_n^T(S)$ and $Y = X^{-1}$. By Lemma~\ref{lem:CommJClass}, $X_{ii} \in U(T)$ for all $i$, as $X_{ii}Y_{ii} = 1_S$. Then, for all $i < j$,
\[ Y_{ii}^{-1}(YX)_{ij} = Y_{ii}^{-1}\sum_{1 \leq k \leq n}Y_{ik}X_{kj}= X_{ij} + Y_{ii}^{-1}\sum_{k \neq i} Y_{ik}X_{kj} = 0_S\]
and hence, $X_{ij} \in V(S)$ as $Y_{ii} \in U(T)$.

\par Now, suppose $X \in UT_n^T(S)$ with $X_{ii} \in U(T)$ and $X_{ij} \in V(S)$ for $i < j$. By \cite[Theorems~3.2 and 4.2]{InvCond}, $X \in UT_n(S)$ is invertible in $M_n(S)$ if and only if $X_{11}^2\cdots X_{nn}^2 \in U(S)$ and $\sum_{k = 1}^n X_{ki}X_{kj} \in V(S)$ for all $i < j$. 
Clearly, $X_{11}^2\cdots X_{nn}^2 \in U(S)$ as $X_{ii} \in U(T)$ and $\sum_{k = 1}^n X_{ki}X_{kj} \in V(S)$ as $X_{ij} \in V(S)$ for all $i \neq j$ and $zx,xz \in V(S)$ for any $x \in V(S)$ and $z \in S$. Thus, $X$ is invertible in $M_n(S)$.

\par Let $Y = X^{-1}$, we aim to show that $Y \in UT_n^T(S)$. Suppose $Y \notin UT_n(S)$ and let $1 < i \leq n$ be the maximum such that there exists $j < i$ with $Y_{ij} \neq 0_S$. Then,
\[ X_{ii}Y_{ij} = (XY)_{ij} = (I_n)_{ij} = 0_S \]
where the first equality holds as $X_{ik} = 0_S$ for all $k < i$ and $Y_{kj} = 0_S$ for all $k > i$ by the maximality of $i$. 
Then, as $X_{ii} \in U(T)$, we get that $Y_{ij} = 0_S$, giving a contradiction, so $Y \in UT_n(S)$. Finally, $Y \in UT_n^T(S)$ as $(XY)_{ii} = X_{ii}Y_{ii} = 1_S$, so $Y_{ii} \in U(T)$.
\end{proof}

\begin{theorem} \label{thm:UTSR}
Let $n \in \N$, $S$ be a commutative semiring, and $T$ be a submonoid of $(S,\cdot)$. Let $\cX$ be a semigroup generating set for the group of units of $UT_n^T(S)$ and let $\Omega,\Theta \subseteq S$ such that $U(T)(\Omega \cup V(S))$ generates $(S,+)$ and $\Theta \cup U(T)$ generates $(T, \cdot)$. Then, the monoid $UT_n^T(S)$ is generated by $\cX \cup E(\Omega) \cup A(\Theta)$ where
\begin{align*}
    A(\Theta) &= \{ A_i(\theta) \colon \theta \in \Theta, \ 1 \leq i \leq n\} \text{, and} \\
    E(\Omega) &= \{ E_{ij}(\omega) \colon \omega \in \Omega, \ 1 \leq i < j \leq n\}.
\end{align*}
Moreover, if $\cX$, $\Omega$ and $\Theta$ are minimal (resp. irredundant) then $UT_n^T(S)$ is minimally (resp. irredundantly) generated by $\cX \cup E(\Omega) \cup A(\Theta)$.
\end{theorem}
\begin{proof}
Let $1 \leq i \leq n$ and $a \in U(T)$, then $A_i(a) \in \lanran{\cX}$ as $A_i(a)$ is invertible by Lemma~\ref{InvCond}.
If $a \in T$, then $a = x_1\cdots x_s$ for some $x_1,\dots,x_s \in \Theta \cup U(T)$. Thus, $A_i(a) =  A_i(x_1)\cdots A_i(x_s)$ and hence $A_i(a)$ is generated by matrices from $A(\Theta) \cup \cX$ for all $1 \leq i \leq n$ and $a \in T$. 

\par Fix $a \in S$. Since $U(T)(\Omega \cup V(S))$ generates $(S,+)$ we can write $a = \sum_{t = 1}^m u_tb_t$ where $u_t \in U(T)$ and $b_t \in \Omega \cup V(S)$. 
Then, for all $i < j$, it is straightforward to verify that 
\[E_{ij}(a) = \prod_{t=1}^{m} A_i(u_t)E_{ij}(b_t)A_i(u_t^{-1}).\]
Moreover, if $b_t \in V(S)$ then $E_{ij}(b_t) \in \lanran{\cX}$ by Lemma~\ref{InvCond}, so $E_{ij}(b_t) \in E(\Omega) \cup \lanran{\cX}$ for all $t$. 
Hence, $E_{ij}(a)$ is generated by $E(\Omega) \cup \cX$ for all $a \in S$ and $i < j$, as $u_t \in U(T)$.
Now, note that, for any $M = (m_{ij}) \in UT_n^T(S)$,
\[ M = \prod_{l=0}^{n-1} \left(A_{n-l}(m_{n-l,n-l})\prod_{k=l+1}^{n-1} E_{n-k,n-l}(m_{n-k,n-l})\right). \]
Therefore, $UT_n^T(S)$ is generated by $\cX \cup E(\Omega) \cup A(\Theta)$.

\par Assume $\cX, \Omega$ and $\Theta$ are minimal and let $\Gamma$ be a generating set for $UT_n^T(S)$ such that $|\Gamma| \leq |\cX \cup E(\Omega) \cup A(\Theta)|$. 
Let $\Gamma_1 \subseteq \Gamma$ be the set of all units in $\Gamma$. 
By Lemma~\ref{Jclass}, any product containing a non-unit is a non-unit, so $\Gamma_1$ generates the group of units, and hence $|\cX| \leq |\Gamma_1|$ as $\cX$ is a minimal generating set for the group of units. Thus, $|\Gamma \setminus \Gamma_1| \leq |E(\Omega) \cup A(\Theta)|$.

\par Let $\cS = \lanran{\cX \cup E(\Omega)}$ and $\Gamma_2 = (\Gamma \setminus \Gamma_1) \cap \cS$. Note that $\cS = UT_n^{U(T)}(S)$ so, $XY \in \cS$ if and only if $X \in \cS$ and $Y \in \cS$ by Lemma~\ref{lem:CommJClass} and hence, $\lanran{\Gamma_1 \cup \Gamma_2} = \cS$. 
We now show that to generate $E_{ij}(x)$ for all $x \in S \setminus V(S)$ and $i < j$, we need at least $|E(\Omega)|$ elements not in $\Gamma_1$.
\par Suppose $\prod_{t=1}^m N_t = E_{ij}(x)$ for some $x \in S \setminus V(S), \ i < j$, and $N_1,\dots,N_m \in UT_n^T(S)$. 
It follows from Lemma~\ref{lem:CommJClass} that $(N_t)_{hh} \in U(T)$ for all $t$ and $h$, since $\prod_{t=1}^m(N_t)_{hh} = (\prod_{t=1}^m N_t)_{hh} = (E_{ij}(x))_{hh} =  1_S$.
So, let $k < l$ such that $(k,l) \neq (i,j)$ then, 
\[\left(\prod_{t=1}^m N_t\right)_{kl} = \sum_{k = i_0 \leq \dots \leq i_m = l} \prod_{s=1}^m (N_s)_{i_{s-1},i_{s}} = (E_{ij}(x))_{kl} = 0_S\]
where the sum is over all possible choices for $i_1, \dots, i_{m-1}$. Thus, for all $1 \leq t \leq m$,
\[(N_1)_{kk}\cdots(N_{t-1})_{kk}(N_t)_{kl}(N_{t+1})_{ll}\cdots(N_m)_{ll} \in V(S) \]
and hence, $(N_t)_{kl} \in V(S)$ as $(N_t)_{hh} \in U(T)$ for all $h$.
Now, for the $(i,j)$ entry, we get that
\[ \left(\prod_{t=1}^m N_t\right)_{ij} =  \sum_{i = i_0 \leq \dots \leq i_m = j} \prod_{s=1}^m (N_s)_{i_{s-1},i_{s}} = (E_{ij}(x))_{ij} = x \]
By the previous paragraph, $(N_t)_{i_l,i_{l+1}} \in V(S)$ if $i_l < i_{l+1}$ and $(i_l,i_{l+1}) \neq (i,j)$. So, we can split this sum in products that contain an entry from $V(S)$ and those that do not.
Let $t_1,\dots,t_{m'}$ be all the indices such that $(N_{t_\alpha})_{ij} \in S \setminus V(S)$ when $1 \leq \alpha \leq m'$, then
\[x = v + \sum_{\alpha=1}^{m'} g_{t_{\alpha}}(N_{t_\alpha})_{ij}\]
for some $v \in V(S)$ and $g_{t_\alpha} \in U(T)$, since the diagonal entries of all $N_t$ are units.

\par Therefore, to generate $E_{ij}(x)$ for all $x \in S \setminus V(S)$, it is necessary to find a set $X \subseteq S$ such that for all $x \in S$, there exist $v \in V(S), \ g_1,\dots,g_{m_x} \in U(T)$, and $x_1,\dots,x_{m_x} \in X$ for some $m_x \in \N_0$ such that $x = v + \sum_{t=1}^{m_x} g_tx_t$.

\par Thus, $U(T)X \cup V(S)$ generates $(S,+)$ and hence, by the definition of $\Omega$, $|X| \geq |\Omega|$, as $U(T)(X \cup V(S)) = U(T)X \cup V(S)$.
Moreover, as we have to generate $E_{ij}(x)$ for all $x \in S \setminus V(S)$ and $i < j$, we get that $|\Gamma_2| \geq \frac{n(n-1)}{2}\cdot |\Omega| = |E(\Omega)|$, and hence $|\Gamma_3| \leq |A(\Theta)|$, where $\Gamma_3 = \Gamma \setminus (\Gamma_1 \cup \Gamma_2)$.

\par For each $s \in T \setminus U(T)$ and $1 \leq i \leq n$, $A_i(s) \notin \lanran{\Gamma_1 \cup \Gamma_2}$, so consider a product $\prod_{t=1}^m N_t = A_i(s)$. Then,
\[(\prod_{t=1}^m N_t)_{ii} = \prod_{t=1}^m (N_t)_{ii} = s \text{ and } (\prod_{t=1}^m N_t)_{hh} = \prod_{t=1}^m (N_t)_{hh} = 1_S\]
for all $h \neq i$. Thus, $(N_t)_{hh} \in U(T)$ for all $t$ and $h \neq i$. Therefore, to generate each $A_i(s)$ for $s \in T \setminus U(T)$ we need to find a set $\Lambda$ such that, for all $s$, there exist $\lambda_1,\dots,\lambda_{m_s} \in \Lambda$ such that $s = g\lambda_1\cdots \lambda_{m_s}$ for some $g \in U(T)$.

\par However, $\Theta$ is the minimal set such that $\Theta \cup U(T)$ generates $(T,\cdot)$, so $|\Lambda| \geq |\Theta|$. Moreover, as we need to generate $A_i(s)$ for all $s \in T \setminus U(T)$ and $1 \leq i \leq n$, we get that $|\Gamma_3| \geq  n|\Theta| = |A(\Theta)|$, and hence $|\Gamma_3| = |A(\Theta)|$. Thus, $|\Gamma| = |\cX \cup E(\Omega) \cup A(\Theta)|$ and $\cX \cup E(\Omega) \cup A(\Theta)$ minimally generates $UT_n^T(S)$.

\par Now, assume $\cX, \Omega$ and $\Theta$ are irredundant. By Lemma~\ref{Jclass}, in $UT_n^T(S)$, any product containing a non-unit is a non-unit. Thus, each element of $\cX$ is irredundant in $\cX \cup E(\Omega) \cup A(\Theta)$.

\par Suppose for a contradiction, $E_{ij}(\omega)$ is redundant for some $i < j$ and $\omega \in \Omega$. Then, to generate $E_{ij}(\omega)$, there exists $v \in V(S), \ g_1,\dots,g_{m_\omega} \in U(T)$, and $x_1,\dots,x_{m_\omega} \in \Omega \setminus \{\omega\}$ for some $m_\omega \in \N_0$ such that $\omega = v + \sum_{t=1}^{m_\omega} g_tx_t$ by above.
This gives a contradiction as $\Omega$ is an irredundant set such that $U(T)(\Omega \cup V(S))$ generates $(S,+)$. 

\par Now, suppose that $A_i(\theta)$ is redundant for some $\theta \in \Theta$.
Then, to generate $A_i(\theta)$, there exist $g \in U(T)$ and $\lambda_1,\dots,\lambda_{m_\theta} \in \Theta \setminus \{\theta\}$ such that $s = g\lambda_1\dots \lambda_{m_\theta}$ by above. This gives a contradiction as $\Theta$ is an irredundant set such that $\Theta \cup U(T)$ generates $(T,\cdot)$.
Thus, $\cX \cup E(\Omega) \cup A(\Theta)$ is a irredundant generating set for $UT_n^T(S)$.
\end{proof}

\begin{remark} 
Let $\cX$, $E(\Omega)$ and $A(\Theta)$ be as defined in the above theorem, so they generate $UT_n^T(S)$. 
Then, $\cX \cup E(\Omega)$ generates $UT_n^{U(T)}(S)$ and $\cX \cup A(\Theta)$ generates $UT^T(V(S))$, that is, the submonoid of $UT^T(S)$ with off-diagonal entries from $V(S)$.
Thus, if $S$ is a ring, then we make take $\Omega$, and hence $E(\Omega)$, to be empty, and if $T = U(T)$, then we may take $\Theta$, and hence $A(\Theta)$, to be empty.
\end{remark}

If we restrict the above theorem to the monoid of $n \times n$ unitriangular matrices over a commutative semiring, we obtain the following corollary describing the generating sets.

\begin{corollary} \label{USR}
Let $n \in \N$ and $S$ be a commutative semiring. Let $\cX$ be a semigroup generating set for the group of units of $UT_n^{\{1_S
\}}(S)$ and let $\Omega \subseteq S$ such that $\Omega \cup V(S)$ generates $(S,+)$. Then, the monoid $UT_n^{\{1_S\}}(S)$ is generated by $\cX \cup E(\Omega)$ where
\[ E(\Omega) = \{ E_{ij}(\omega) \colon \omega \in \Omega, \ 1 \leq i < j \leq n\}. \]
Moreover, if $\cX$ and $\Omega$ are minimal (resp. irredundant) then $UT_n^{\{1_S\}}(S)$ is minimally (resp. irredundantly) generated by $\cX \cup E(\Omega)$.
\end{corollary}

If we apply the above two results in the case when $S = \Zmax$, we obtain the following corollaries constructing explicit minimal and irredundant generating sets for $UT_n(\Zmax)$ and $UT_n^{\{0\}}(\Zmax)$.

\begin{corollary} \label{UTNZ}
Let $n \in \N$. Then, the monoid $UT_n(\Zmax)$ is minimally and irredundantly generated by $A(1) \cup \{-1 \cdot I_n\} \cup E(0) \cup A(-\infty)$ where
\[
A(1) = \{A_i(1) \colon 1 \leq i \leq n \}, \ E(0) = \{E_{ij} \colon 1 \leq i <  j \leq n \} \text{, and }\]
\[ A(-\infty) = \{A_i(-\infty) \colon 1 \leq i \leq n \}. \]
\end{corollary}
Recall that $1 \neq 1_{\Zmax} = 0 \neq 0_{\Zmax} = -\infty$ and that $-1 \cdot I_n$ is the diagonal matrix with $-1$ on the diagonal and $-\infty$ elsewhere.
\begin{proof}
As $\Zmax$ is an anti-negative semifield, $U(\Zmax) = \Z$ and $V(\Zmax) = \{-\infty\}$, so, by Lemma~\ref{InvCond}, $X \in UT_n(\Zmax)$ is invertible if and only if $X$ is diagonal with $X_{ii} \neq -\infty$ for all $i$. Moreover, $\Z(\{0\} \cup \{-\infty\}) = \Zmax$ and $\{-\infty\} \cup \Z = \Zmax$ so, by Theorem~\ref{thm:UTSR}, it suffices to show that $\cX = \{-1 \cdot I_n, A_1(1), \dots A_n(1)\}$ forms a minimal and irredundant generating set for the group of units of $UT_n(\Zmax)$.

Clearly, the group of units of $UT_n(\Zmax)$ is isomorphic to $\Z^n$ under coordinate-wise addition and is generated by $\cX$. Finally, we can see that $\cX$ is minimal and irredundant as $|\cX| = n+1$ and $\Z^n$ is minimally $(n+1)$-generated as a semigroup \cite[Corollary 4.3]{SemRank}.
\end{proof}

\begin{corollary}
Let $n \in \N$. Then, the monoid $UT_n^{\{0\}}(\Zmax)$ is minimally and irredundantly generated by $\{I_n\} \cup E(\Z)$ where
\[ E(\Z) = \{E_{ij}(z) \colon z \in \Z, 1 \leq i <  j \leq n \}. \]
\end{corollary}
\begin{proof}
Remark that $\max(x,y) \in \{x,y\}$ for all $x,y \in \Zmax$. Thus, the minimal and irredundant generating set for $(\Z,\max)$ is $\Z$.
\end{proof}

\section{Generating sets for full matrix monoids} \label{FullGenSec}
We now focus on constructing generating sets of full matrix monoids over anti-negative semifields. In particular, we provide minimal and irredundant generating sets for $M_2(\Zmax)$ and $M_3(\Zmax)$, showing that $M_n(\Zmax)$ is finitely generated if and only if $n \leq 2$.

\par We define two functions which we use throughout this section. For a semiring $S$, define $\psi \colon S \rightarrow \B$ to be the map that sends $0_S$ to $0$ and $S^*$ to $1$, and $\phi_n \colon M_n(S) \rightarrow M_n(\B)$ to be the map where $\phi_n(A)_{ij} = \psi(A_{ij})$. 
If $S$ is a non-trivial anti-negative semiring without zero-divisors, then $\psi$ and $\phi_n$ are surjective morphisms for all $n \in \N$, and hence the cardinality of a minimal generating set for $M_n(S)$ is at least the cardinality of a minimal generating set for $M_n(\B)$.
\subsection{{2-by-2} full matrix monoids}

\par We say $M \in M_n(S)$ is a \emph{monomial matrix} if there exists $\sigma \in \cS_n$ such that $M_{ij} \neq 0_S$ if and only if $j = \sigma(i)$, and we say that $M$ has \emph{underlying permutation} $\sigma$. Moreover, a monomial matrix $M$ is the \emph{permutation matrix of} $\sigma$ if $M_{ij} = 1_S$ for all $j = \sigma(i)$. 

The following lemma tells us when a matrix over a commutative anti-negative semiring without zero divisors is invertible, this can be deduced from \cite[Corollary 3.3]{AntiNegInvCond}. We denote the group of units of $M_n(S)$ as $GL_n(S)$.
\begin{lemma} \label{Invert}
Let $S$ be a commutative anti-negative semiring without zero divisors. Then, $GL_n(S)$ consists exactly of the monomial matrices where all non $0_S$ entries are in $U(S)$.
\end{lemma}
\begin{proof}
Note that monomial matrices in which every non $0_S$ entry is in $U(S)$ satisfy the conditions of \cite[Corollary 3.3]{AntiNegInvCond} and hence are invertible. So, now suppose that $X \in M_n(S)$ is invertible. Then, by \cite[Corollary 3.3(2--3)]{AntiNegInvCond}, we can see that $X_{ij}X_{ik} = 0_S = X_{ji}X_{ki}$ for all $1 \leq i,j,k \leq n$ with $j \neq k$. Thus, as $S$ has no zero-divisors, $X$ has at most one non $0_S$ entry per row and column. Finally, observe that, by \cite[Corollary 3.3(2)]{AntiNegInvCond}, all non $0_S$ entries of $X$ are in $U(S)$.
\end{proof}

For a semiring $S$, we say that, for $x,y \in S$, $x \leq y$ if and only if there exists $t \in S$ such that $x + t = y$. We say $S$ is \emph{linearly ordered} if $x \leq y$ or $y \leq x$ for all $x,y \in S$.

\begin{theorem} \label{thm:2by2SemiField}
Let $S$ be a linearly ordered anti-negative semifield. Let $X$ be a semigroup generating set for $(S^*,\cdot)$. If $(S^*,\cdot)$ is non-trivial, choose $X$ such that $\alpha^{-1} \in \lanran{X \setminus \{\alpha\}}$ for some $\alpha \in X$.
Then, the monoid $M_2(S)$ is generated by the matrices:
\[A_1(x) =
\begin{pmatrix}
x & 0_S \\
0_S & 1_S
\end{pmatrix} \text{ for all } x \in X \setminus \{\alpha\},
\]
\[ 
B =
\begin{pmatrix}
0_S & \alpha \\
1_S & 0_S
\end{pmatrix}, \
C = 
\begin{pmatrix}
0_S & 0_S \\
0_S & 1_S
\end{pmatrix}
\text{, and }
D =
\begin{pmatrix}
1_S& 1_S\\
0_S& 1_S
\end{pmatrix}\]
Moreover, if $X$ is minimal (resp. irredundant), then $M_2(S)$ is minimally (resp. irredundantly) generated by $A_1(x)$, $B$, $C$, and $D$ for $x \in X \setminus \{\alpha\}$.
\end{theorem}
\begin{proof}
We begin by noting that, by our choice of $X$, when $(S^*,\cdot)$ is non-trivial, there exists $x_1,\dots,x_s \in (X \setminus \{\alpha\})$ such that $x_1\cdots x_s = \alpha^{-1}$. Thus, $A_1(\alpha^{-1}) = A_1(x_1)\cdots A_1(x_s)$ and when $(S^*,\cdot)$ is trivial, $B^2 = A_1(\alpha^{-1}) = I_2$. Thus, in either case, we can generate $A_1(\alpha^{-1})$, and hence also,
\[
F = 
\begin{pmatrix}
0_S & 1_S\\
1_S& 0_S
\end{pmatrix} = A_1(\alpha^{-1})B \text{ and } 
A_1(\alpha) =
\begin{pmatrix}
\alpha & 0_S \\
0_S & 1_S
\end{pmatrix} = BA_1(\alpha^{-1})B.
\]
Thus, we can generate $A_1(z)$ for all $z \in S^*$, as $X$ generates $(S^*,\cdot)$, so $A_1(z) =  A_1(x_1)\cdots A_1(x_t)$ for some  $x_1,\dots,x_t \in X$.
Moreover, pre-multiplying a matrix by $F$ swaps the rows and post-multiplying by $F$ swaps the columns, so it suffices to show that we can generate every matrix, up to rearranging rows and columns. 
Now, observe that, for $x,y,z \in S^*$,
\begin{align*}
\begin{pmatrix}
0_S & 0_S \\
0_S & 0_S
\end{pmatrix} &= CFC,
&&\begin{pmatrix}
0_S & 0_S \\
x & 0_S
\end{pmatrix} = CFA_1(x), \\
\begin{pmatrix}
0_S & 0_S \\
x & y 
\end{pmatrix} &= CFDA_1(y)FA_1(x), 
&&\begin{pmatrix}
0_S & x \\
0_S & y
\end{pmatrix} = A_1(x)FA_1(y)DC, \\
\begin{pmatrix}
0_S & x \\
y & 0_S
\end{pmatrix} &= A_1(x)FA_1(y) \text{, and} 
&&\begin{pmatrix}
0_S & x \\
y & z
\end{pmatrix} = A_1(x)FA_1(z)DA_1(z^{-1}y).
\end{align*}
Therefore, every matrix with at least one $0_S$ entry is a product of the given matrices. Finally, for $a,b,c,d \in S^*$, note that
\[
\begin{pmatrix}
a & b \\
c & d
\end{pmatrix} =
\begin{pmatrix}
1_S& 1_S\\
db^{-1} & ca^{-1}
\end{pmatrix}
A_1(b)FA_1(a). \]
So, it suffices to express
$\left(\begin{smallmatrix}
1_S& 1_S\\
x & y 
\end{smallmatrix}\right)$
as a product of matrices with at least one $0_S$ entry for all $x,y \in S^*$. 
Without loss of generality, we may suppose $y \leq x$ as if $x \leq y$, then we can post-multiply by $F$ to swap the columns.
So, as $y \leq x$, there exists $t \in S$ such that $t + y = x$ and 
\[
\begin{pmatrix}
1_S& 1_S\\
x & y 
\end{pmatrix} =
\begin{pmatrix}
0_S & 1_S \\
y & y
\end{pmatrix}
\begin{pmatrix}
y^{-1}t & 0_S \\
1_S& 1_S
\end{pmatrix}.\]
Thus, every matrix without $0_S$ entries is a product of matrices with at least one $0_S$ entry and hence, $M_2(S)$ is generated by the given matrices.
\par Now, we show that if $X$ is minimal then the generating set is minimal. 
By Lemma~\ref{Invert}, $GL_2(S)$ is the set of monomial matrices with entries in $S^*$. 
So, let $\mathrm{perm} \colon GL_2(S) \rightarrow (S^*, \cdot)$ be the surjective morphism that maps a matrix to the product of its non $0_S$ entries.

\par As $(S^*, \cdot)$ is minimally generated by $X$, $GL_2(S)$ is minimally generated by at least $|X|$ matrices. Moreover, any generating set for $M_2(S)$ contains a generating set for $GL_2(S)$ by Lemma~\ref{Jclass}. Thus, if $X$ is infinite then we are done, so assume $X$ is finite.

\par Now, for a contradiction, suppose there exists a generating set $\Gamma$ of size $|X| + 1$ for $M_2(S)$. By above $|X|$ elements of $\Gamma$ are in $GL_2(S)$. 
Let $\Gamma' = \Gamma \cap GL_2(S)$ and $\Gamma \setminus \Gamma' = \{\gamma\}$. 
Moreover, as $\phi_2$ is a surjective morphism, $\phi_2(\Gamma')$ and $\phi_2(B)$ are generating sets for $GL_2(\B)$, and $\phi_2(\Gamma)$ is a generating set for $M_2(\B)$.
Thus, $\phi_2(B) \cup \phi_2(\gamma)$ is a generating set for $M_2(\B)$, giving a contradiction as $M_2(\B)$ is minimally generated by 3 matrices \cite[Table 1]{MinGens}. Therefore, $M_2(S)$ is minimally generated by these $|X|+2$ matrices. 

\par Now, suppose $X$ is irredundant. By again considering the surjective morphism $\mathrm{perm} \colon GL_2(S) \rightarrow (S^*,\cdot)$, we see that $B$ and $A_1(x)$ for all $x \in X \setminus \{\alpha\}$ form an irredundant generating set for $GL_2(S)$. 
By Lemma~\ref{Jclass}, any generating set for $M_2(S)$ contains a generating set for $GL_2(S)$, so $B$ and $A_1(x)$ for all $x \in X \setminus \{\alpha\}$ are also irredundant in the generating set for $M_2(S)$.
Finally, $\phi_2(A_1(x)) = I_2$ for all $x \in X \setminus \{\alpha\}$ and $M_2(\B)$ is minimally, and hence irredundantly, generated by $\phi_2(B)$, $\phi_2(C)$, and $\phi_2(D)$, \cite[Table 1]{MinGens}, so $C$ and $D$ are irredundant. 
\end{proof}

The following proposition shows that, by using the above theorem, we can always find a minimal generating set for $M_2(S)$ when $S$ is a linearly ordered anti-negative semifield.

\begin{proposition}
    Let $S$ be a linearly ordered anti-negative semifield with $(S^*,\cdot)$ non-trivial. Then, there exists a minimal semigroup generating set $X$ for $(S^*,\cdot)$ such that $\alpha^{-1} \in \lanran{X \setminus \{\alpha\}}$ for some $\alpha \in X$.
\end{proposition}
\begin{proof}
    Clear, if $(S^*,\cdot)$ is not finitely generated. So suppose $(S^*,\cdot)$ is finitely generated. As $S$ is an anti-negative semifield, every element but $0_S$ and $1_S$ has infinite multiplicative order \cite[Lemma 2.1(ii)]{MatrixOverSemirings}, and hence, $(S^*,\cdot)$ is isomorphic to $\Z^m$ for some $m \in \N$. 
So, let $X' = \{x_1, \dots, x_m\}$ be a minimal group generating set for $(S^*,\cdot)$ and $X = X' \cup \{x_1^{-1}\cdots x_m^{-1}\}$. 
Then, $X$ is a minimal semigroup generating set such that $x_1^{-1} \in \lanran{X \setminus \{x_1\}}$ as $\Z^m$ is minimally generated by $m+1$ elements as a semigroup \cite[Corollary 4.3]{SemRank}.

\end{proof}

\begin{corollary} \label{2by2Z}
The monoid $M_2(\Zmax)$ is minimally generated by:
\[
A =
\begin{pmatrix}
1 & -\infty \\
-\infty & 0
\end{pmatrix}, \
B =
\begin{pmatrix}
-\infty & -1 \\
0 & -\infty
\end{pmatrix},
\]
\[ 
C = 
\begin{pmatrix}
-\infty & -\infty \\
-\infty & 0
\end{pmatrix}
\text{, and }
D =
\begin{pmatrix}
0 & 0 \\
-\infty & 0
\end{pmatrix}\]
\end{corollary}
\begin{proof}
Note that $\Zmax$ is linearly ordered, $X = \{-1,1\}$ is a generating set for $(\Z,+)$, and $(-1)^{-1} = 1$.
\end{proof}

\subsection{Higher dimension full matrix monoids}
We now turn our attention to the monoids $M_n(S)$ where $n \geq 3$. In particular, we show that there are infinitely many prime $\cJ$-classes in $M_n(S)$ when $n \geq 3$ and $S$ is an infinite commutative anti-negative semiring without zero divisors, and hence that $M_n(S)$ is not finitely generated.
We then construct an (infinite) minimal and irredundant generating set for $M_3(\Zmax)$.
\par First, we introduce notation for a collection of matrices which we use for the remainder of this section. For a semiring $S$, $n \geq 2$, and $s \in S^*$, let $Z_n(s) \in M_n(S)$ be 
\[ 
Z_n(s) =
\begin{pmatrix}
1_S & 1_S & 0_S & \cdots & 0_S \\
0_S & \ddots & \ddots & \ddots & \vdots \\
\vdots & \ddots & \ddots & \ddots & 0_S \\
0_S & \ddots & \ddots & \ddots & 1_S \\
s & 0_S & \cdots & 0_S & 1_S \\
\end{pmatrix},
\]
that is, $Z_n(s)_{ij} = 1_S$ if $j = i$ or $j = i+1$, $Z_n(s)_{n1} = s$, and $Z_n(s)_{ij} = 0_S$ otherwise.

\begin{lemma} \label{primeJ}
Let $n \geq 3$ and $S$ be a commutative anti-negative semiring without zero divisors. 
Then, $Z_n(s)$ is prime in $M_n(S)$ for all $s \in S^*$. Moreover, if $Z_n(s) \cJ Z_n(t)$ for some $t \in S^*$ then $s = t$ or $st = 1_S$. 
\end{lemma}
\begin{proof}
By \cite[Theorem 1]{PrimeBoolean}, $\phi_n(Z_n(s))$ is prime, so, $AB = Z_n(s)$ implies either $\phi_n(A)$ or $\phi_n(B) $ is a unit. 
If $\phi_n(A)$ is a unit, then it is a permutation matrix and $A$ is a monomial matrix.
Hence, if $A$ is not a unit then, $A$ has a non $0_S$, non-invertible entry by Lemma~\ref{Invert}. 
Thus, some row of $AB$ is a scaling of a row of $B$ by a non-invertible element of $S$. However, $1_S$ is an entry of each row of $Z_n(s)$, giving a contradiction by Lemma~\ref{lem:CommJClass}. Hence, $A$ is a unit. 
If $\phi_n(B)$ is a unit a dual argument holds, since $1_S$ is an entry of each column of $Z_n(s)$. Therefore, $Z_n(s)$ is prime in $M_3(S)$.

\par Let $\cX = \{ X \in M_n(S) \colon \phi_n(X) = \phi_n(Z_n(1_S))\}$ and define $v \colon \cX \rightarrow S^* \times S^*$ to be the map, where $v(X) = ((X_{1,2}\cdots X_{n-1,n}X_{n,1}), (X_{1,1}\cdots X_{n,n}))$. 
Note that $v(Z_n(s)) = (s,1_S)$ for all $s \in S^*$. Say $X \equiv Y$ for $X,Y \in \cX$ if there exists $g \in U(S)$ such that $v(X) = (g,g) \cdot v(Y)$.

\par Let $s,t \in S^*$ and suppose $Z_n(s) \cJ Z_n(t)$. Then, there exists $U,V \in GL_n(S)$ such that $UZ_n(s)V = Z_n(t)$, as $Z_n(t)$ is prime. 
By Lemma~\ref{Invert}, we may write $U = DP$ and $V = P'D'$ for permutation matrices $P$ and $P'$ and diagonal matrices with entries in $U(S)$, $D$ and $D'$. 
\par Let $X \in \cX$, and consider $DXD'$. Each entry of $D$ and $D'$ scales a row and column of $X$ respectively, and hence, scale one entry from both $\{X_{1,2},\dots, X_{n-1,n},X_{n,1}\}$ and $\{X_{1,1},\dots,X_{n,n}\}$ by some $d \in U(S)$.
Thus, $v(DXD') = (g,g) \cdot v(X)$ for some $g \in U(S)$. Therefore, $X \equiv DXD'$ and $PZ_n(s)P' \equiv Z_n(t)$. Moreover, as $\phi_n(PZ_n(s)P') = \phi_n(Z_n(t))$, we only have the consider the permutations of $Z_n(s)$ contained in $\cX$.

\par Let $Y = PZ_n(s)P'$, then for some $1 \leq i \leq n$, either $Y_{i,i} = s$, $Y_{i,i+1} = s$ or $Y_{n,1} = s$ with all other non $0_S$ entries equal to $1_S$. 
Hence, $v(PZ_n(s)P') = (s,1_S)$ or $(1_S,s)$. 
Therefore, $PZ_n(s)P' \equiv Z_n(t)$ implies that $s = t$ or $st = 1_S$, and hence, if $Z_n(t) \cJ Z_n(s)$ then $s = t$ or $st = 1_S$.
\end{proof}

\begin{theorem} \label{InfConj}
Let $n \geq 3$ and $S$ be an infinite commutative anti-negative semiring without zero divisors. Then, the monoid $M_n(S)$ is not finitely generated.
\end{theorem}
\begin{proof}
Let $\cZ = \{ Z_n(s) \colon s \in S^*\}$. By Lemma~\ref{primeJ}, each $Z \in \cZ$ is prime and hence any generating set for $M_n(S)$ contains a matrix $\cJ$-related to each $Z \in \cZ$. 
However, each $Z_n(s) \in \cZ$, is $\cJ$-related to at most one other matrix from $\cZ$ by Lemma~\ref{lem:CommJClass} and \ref{primeJ}. Thus, $M_n(S)$ is not finitely generated, as $S$, and hence $\cZ$ is infinite.
\end{proof}

By Theorem~\ref{thm:2by2SemiField} and \ref{InfConj} we obtain the following corollary.
\begin{corollary} \label{M3ZNotFin} 
The monoid $M_n(\Zmax)$ is finitely generated if and only if $n \leq 2$.
\end{corollary}

For the remainder of this section, when the dimension of the matrix is clear, we use the notation $P_\sigma \in M_n(S)$ for the permutation matrix of $\sigma \in \cS_n$.
\begin{lemma} \label{units}
Let $S$ be a commutative anti-negative semiring without zero divisors, $X$ be a generating set for $(U(S),\cdot)$, and $x \in X$. Then, for $n \geq 2$, $GL_n(S)$ is generated by
\[ A = A_1(x)P_{(1,\dots,n-1)}, \ B = A_1(x^{-1})P_{(1,\dots,n)}, \text{ and } A_1(y) \text{ for } y \in X \setminus \{x,x^{-1}\}.\]
\end{lemma}
\begin{proof}
\par Remark that $A^{n-1} = A_1(x)\cdots A_{n-1}(x)$. Then,
\begin{align*}
B^{n-2}A^{n-1}B &= (A_1(x^{-1})P_{(1,\dots,n)})^{n-2}A_1(x)\cdots A_{n-1}(x)A_1(x^{-1})P_{(1,\dots,n)} \\ 
&= P_{(1,\dots,n)}^{n-2}A_{n-1}(x^{-1})\cdots A_{2}(x^{-1})A_2(x)\cdots A_{n-1}(x)P_{(1,\dots,n)}  \\
&= P_{(1,\dots,n)}^{-1}
\end{align*}
as $A_i(x^{-1})P_{(1,\dots,n)} = P_{(1,\dots,n)}A_{i+1}(x^{-1})$ for all $i$. Therefore, it follow that, $(B^{n-2}A^{n-1}B)^{n-1} = P_{(1,\dots,n)}$. Moreover,
\begin{align*}
BP_{(1,\dots,n)}^{-1}A = A_1(x^{-1})A_1(x)P_{(1,\dots,n-1)} &= P_{(1,\dots,n-1)} \text{, and} \\
P_{(1,\dots,n)}^{-2}P_{(1,\dots,n-1)}P_{(1,\dots,n)} &= P_{(1,2)}.
\end{align*} 
Thus, every permutation matrix is generated by $A$ and $B$, as $\cS_n$ is generated by the permutations $(1,2)$ and $(1,\dots,n)$ \cite[Exercise 2.9(iii)]{Rotman2012}. Moreover,
\[A_i(x) = P_{(1,i)}AP_{(1,\dots,n-1)}^{-1}P_{(1,i)}, \ A_i(x^{-1}) = P_{(1,i)}BP_{(1,\dots,n)}^{-1}P_{(1,i)}, \]
\[ \text{and } A_i(y) = P_{(1,i)}A_1(y)P_{(1,i)} \text{ for } y \in X \setminus \{x,x^{-1}\}.\]
Hence, every diagonal matrix with entries in $U(S)$ can be generated, as they can be expressed as a product using matrices $A_i(x)$ for $x \in X$ where $1 \leq i \leq n$, which can be seen to be generated by above.

\par By Lemma~\ref{Invert}, every matrix in $GL_n(S)$ can be expressed as diagonal matrix with entries from $U(S)$ multiplied by a permutation matrix. Therefore, $GL_n(S)$ is generated by $A$, $B$, and $A_1(y) \text{ for } y \in X \setminus \{x,x^{-1}\}$.
\end{proof}

\begin{corollary} \label{unitsZ}
Let $n \geq 2$. $GL_n(\Zmax)$ is minimally generated by the matrices $A = A_1(1)P_{(1,\dots,n-1)} \text{ and } B = A_1(-1)P_{(1,\dots,n)}$.
\end{corollary}
\begin{proof}
The group $GL_n(\Zmax)$ is non-abelian and hence not generated by one matrix. Thus, by Lemma~\ref{units}, $A$ and $B$ minimally generate $GL_n(\Zmax)$.
\end{proof}

\begin{lemma} \label{lem:JRelatedtoZ}
    Let $S$ be a semifield and suppose $X \in M_3(S)$ has exactly one $0_S$ in each row and column. Then, $X \cJ Z_3(s)$ for some $s \in S^*$.
\end{lemma}
\begin{proof}
    Let $0 = 0_S$, $1 = 1_S$, and $\cX \subseteq M_3(S)$ contain all matrices with exactly one $0$ in each row and column. Then, by multiplying by permutation matrices $X$ is $\cJ$-related to a matrix $Y$ with $Y_{ii} = 0$ for all $i$. Finally, note that,
    \[ 
\begin{pmatrix}
    0 & 1 & 1 \\
    1 & 0 & 1 \\
    s & 1 & 0
\end{pmatrix} 
= 
\begin{pmatrix}
    b^{-1} & 0 & 0 \\
    0 & d^{-1} & 0 \\
    0 & 0 & ab^{-1}f^{-1} \\
\end{pmatrix} 
\begin{pmatrix}
    0 & a & b \\
    c & 0 & d \\
    e & f & 0
\end{pmatrix} 
\begin{pmatrix}
    c^{-1}d & 0 & 0 \\
    0 & a^{-1}b & 0 \\
    0 & 0 & 1 \\
\end{pmatrix}
\]
where $s = ab^{-1}c^{-1}def^{-1}$. Thus, $X \cJ Y \cJ Z_3(s)$.
\end{proof}

\par For matrices $X,Y \in M_n(S)$, we say that $X$ is a \emph{permutation} of $Y$ if $X$ can be obtained by permuting the rows and columns of $Y$. Equivalently, $X = PYP'$ for some permutation matrices $P,P' \in M_n(S)$.

\begin{theorem} \label{thm:4regularmatrices3by3}
    Let $\cX \subseteq M_3(\Zmax)$ contain all matrices with exactly one $-\infty$ entry in each row and column. Then, the submonoid $M_3(\Zmax) \setminus \cX$ is minimally and irredudantly generated by 
    \[ A = A_1(1)P_{(1,2)}, \ B = A_1(-1)P_{(1,2,3)}, \ E_{12} \text{, and } A_1(-\infty).\]
\end{theorem}
\begin{proof}
Note that $\cX$ only contains prime $\cJ$-classes by Lemma~\ref{primeJ} and \ref{lem:JRelatedtoZ}, so $\cX$ is a submonoid of $M_3(\Zmax)$. Let $\cR$ denote the monoid generated by $A,B,E_{12}$, and $A_1(-\infty)$, so we aim to show $\cR = M_3(\Zmax) \setminus \cX$.

\par By Corollary~\ref{unitsZ}, $A$ and $B$ generated $GL_3(\Zmax)$, so $GL_3(\Zmax) \subseteq \cR$. 
Thus, $-1 \cdot I_3 \in \cR$ and $A_i(1) \in \cR$ for $1 \leq i \leq 3$. Moreover, for $1 \leq k \leq 3$ and $i < j$, observe that
\[ A_k(-\infty) = P_{(k1)}A_1(-\infty)P_{(1k)} \text{, and }  E_{ij} = P_{(1i)(2j)}E_{12}P_{(2j)(1i)}.
\]
Thus, $UT_3(\Zmax) \subseteq \cR$ as $P_\sigma \in GL_3(\Zmax)$ for all $\sigma \in \cS_3$, and hence, all the generators from Corollary~\ref{UTNZ} are contained in $\cR$.

Note that $A_1(1)$, $A_1(-1)P_{(12)}$, $A_1(-\infty)$, and $E_{12}$ are contained in $\cR$. Moreover, by restricting these matrices to their first two rows and columns, we obtain a generating set for $M_2(\Zmax)$ by Lemma~\ref{2by2Z}, as they are block diagonal and the $(3,3)$ entry of each matrix is $0$.
Hence, by multiplying by $A_3(x) \in \cR$ for $x \in \Zmax$, we can construct any block diagonal matrix with a $2 \times 2$ block and a $1 \times 1$ block.

\par If a matrix has at least four $-\infty$ entries then it contains a row and column with at least two $-\infty$ entries. Thus, every matrix with at least four $-\infty$ entries is contained in $\cR$ as it is a permutation of either an upper triangular matrix or a block diagonal matrix with a $2 \times 2$ block.

\par Next, we show that we can construct all matrices with three $-\infty$ entries apart from those in $\cX$. Note that $A^T,B^T, (E_{12})^T$, $A_1(-\infty)^T \in \cR$ as they have more than four $-\infty$ entries, and hence, we only have to show we can generate all matrices up to transposition and permutation.
\par Now, for $a,b,c,d,e,f,g \in \Zmax$,
\[
\begin{pmatrix}
a & b & c \\
d & e & f \\
-\infty & -\infty & g
\end{pmatrix} =
\begin{pmatrix}
0 & -\infty & c\\
-\infty & 0 & f \\
-\infty & -\infty & g
\end{pmatrix}
\begin{pmatrix}
a & b & -\infty \\
d & e & -\infty \\
-\infty & -\infty & 0
\end{pmatrix}.
\]
Thus, the above matrix is the product of matrices with at least four $-\infty$ entries, so $\cR$ contains any matrix with at least two $-\infty$ entries in the same row or column. Thus, $\cR$ contains all matrices not in $\cX$ with at least three $-\infty$ entries.
\par Now, for $a,b,c,d,e,f,g \in \Z$ and $x \in \Zmax$, 
\[
\begin{pmatrix}
a & b & c \\
d & e & -\infty \\
f & x & g
\end{pmatrix} =\begin{cases}
\left(\begin{smallmatrix}
0 & b-e & -\infty \\
-\infty & 0 & -\infty \\
-\infty & -\infty & 0
\end{smallmatrix}\right)
\left(\begin{smallmatrix}
a & -\infty & c \\
d & e & -\infty \\
f & x & g
\end{smallmatrix}\right), &\text{if } a + e \geq b + d \vspace{2mm} \\ 
\left(\begin{smallmatrix}
a-d & 0 & -\infty \\
0 & -\infty & -\infty \\
-\infty & -\infty & 0
\end{smallmatrix}\right)
\left(\begin{smallmatrix}
d & e & -\infty \\
-\infty & b & c \\
f & x & g
\end{smallmatrix}\right), &\text{if } b + d \geq a + e.
\end{cases} 
\]
By taking $x = -\infty$ above, we can see that $\cR$ contains all matrices with two $-\infty$ entries as they are the product of matrices with at least three $-\infty$ entries not in $\cX$. Similarly, taking $x \in \Z$, shows that $\cR$ contains any matrix with one $-\infty$ entry as they product matrices with at least two $-\infty$ entries not in $\cX$.

\par Finally, for matrices without $-\infty$ entries, we may scale the columns so that the top row only contains $0$ entries. So, we only need to consider matrices of the form
\[
\begin{pmatrix}
0 & 0 & 0 \\
a & b & c \\
d & e & f
\end{pmatrix}\] 
where $a,b,c,d,e,f \in \Z$. Further, we may rearrange the columns to assume $a \leq b,c$ and $e \leq f$, then
\[ 
\begin{pmatrix}
0 & 0 & 0 \\
a & b & c \\
d & e & f
\end{pmatrix} = \begin{cases}
\left(\begin{smallmatrix}
0 & -\infty & -\infty \\
a & b & c \\
d & e & f
\end{smallmatrix}\right)
\left(\begin{smallmatrix}
0 & 0 & 0 \\
-\infty & 0 & -\infty \\
-\infty & -\infty & 0
\end{smallmatrix}\right) &\text{if } d \leq e, \vspace{2mm} \\
\left(\begin{smallmatrix}
0 & -b & -d \\
c & 0 & -\infty \\
f & -\infty & 0
\end{smallmatrix}\right)
\left(\begin{smallmatrix}
-\infty & -\infty & 0 \\
a & b & -\infty \\
d & e & -\infty
\end{smallmatrix}\right) &\text{if } e \leq d, 
\end{cases}\]
where the second case holds as $a - b \leq 0$ and $e - d \leq 0$. 
Thus, $\cR = M_3(\Zmax) \setminus \cX$.
To see that the generating set is minimal, note that the $GL_3(\Zmax)$ is non-abelian, and hence requires at least 2 matrices to generate it. Moreover, to generate $E_{12}$ we require a permutation of $E_{12}$, and hence a matrix with a row only containing $-\infty$ is required, as $A,B$ and $E_{12}$ have an entry in each row and column. Thus, the generating set is minimal and hence, irredundant.
\end{proof}
\par For a semigroup $\cS$, we say $x \in \cS$ is \emph{regular} if there exists $y \in \cS$ such that $xyx = x$. 
In 1968, Devadze \cite{Devadze} showed that the size of minimal generating sets for $M_n(\B)$ grows exponentially. However, Kim and Roush \cite{KimRoush} showed that, for all $n \in \N$, there exists a subsemigroup of $M_n(\B)$ generated by four matrices which contains all regular matrices in $M_n(\B)$. 

\par Note that, as $\cX$ from the above theorem only contains prime $\cJ$-classes, it is not hard to show $\cX$ contains no regular matrices. Thus, there exists a subsemigroup of $M_3(\T)$ generated by four matrices which contains all regular matrices in $M_3(\T)$. 
We now pose the question of whether the theorem is true when applied to $M_n(\Zmax)$ for all $n \in \N$.
\begin{question}
Do the matrices 
\[ A_1(1)P_{(1,\dots,n-1)}, \ A_1(-1)P_{(1,\dots,n)}, \ E_{12} \text{, and } A_1(-\infty)\]
generate all regular matrices of $M_n(\Zmax)$?
\end{question}

Finally, by adjoining an element from each prime $\cJ$-class to our generating set above, we obtain an (infinite) minimal and irredundant generating set for $M_3(\Zmax)$.

\begin{corollary} \label{cor:FM3Z}
The monoid $M_3(\Zmax)$ is minimally and irredundantly generated by the following matrices:
\[A = A_1(1)P_{(1,2)}, \ B = A_1(-1)P_{(1,2,3)}, \ E_{12}, \ A_1(-\infty) \text{, and } Z_3(i) \text{ for } i \in \N_0
\]
\end{corollary} 
\begin{proof}
By Theorem~\ref{thm:4regularmatrices3by3}, it suffices to show that we can generate all matrices with exactly one $-\infty$ in each row and column. By Lemma~\ref{lem:JRelatedtoZ}, each of these matrices is $\cJ$-related to $Z_3(i)$ for some $i \in \Z$. Moreover, by Lemma~\ref{primeJ}, each $Z_3(i)$ is prime, so each of the matrices can be obtained by multiplying some $Z_3(i)$ by matrices in $GL_3(\Zmax)$. Thus, it suffices to show we can generate $Z_3(i)$ for all $i \in \Z$.

\par Clearly, we can generate $Z_3(i)$ for $i \in \N_0$, so suppose $i < 0$. Then, 

\[ Z_3(i) = 
\begin{pmatrix}
    0 & -\infty & -\infty \\
    -\infty & -\infty & i \\
    -\infty & i & -\infty 
\end{pmatrix} Z_3(-i)
\begin{pmatrix}
    -\infty & 0 & -\infty \\
    0 & -\infty & -\infty \\
    -\infty & -\infty & -i
\end{pmatrix} 
\]
Thus, each matrix with exactly one $-\infty$ in each row and column can be generated. Hence, the given matrices form a generating set for $M_3(\Zmax)$.

\par The generating set is minimal by Corollary~\ref{M3ZNotFin}. We now show that the generating set is irredundant. Note that each $Z_3(i)$ is irredundant as, by Lemma~\ref{primeJ}, they are all contained in different prime $\cJ$-classes, and every generating set requires a representative from each.

\par Moreover, $\phi_3(A)$, $\phi_3(B)$, $\phi_3(E_{12})$, $\phi_3(A_1(-\infty))$, and $\phi_3(Z_3(0))$ provide a generating set for $M_3(\B)$ as $\phi_3$ is a surjective morphism and $\phi_3(Z_3(i)) = \phi_3(Z_3(0))$ for all $i \in \N_0$. However, $M_3(\B)$ is minimally generated by 5 matrices \cite[Table 1]{MinGens}, so $A$, $B$, $E_{12}$, and $A_1(-\infty)$ are irredundant. Therefore, the generating set is irredundant. 
\end{proof}

\section{Presentations of tropical matrix monoids} \label{PresentationSec}
In this section, we establish that $UT_n(\Zmax)$ is finitely presented for all $n \in \N$, and we construct an explicit presentation using the minimal generating set from Section~\ref{UTGenSec}. We then show that $M_n(\Zmax)$ is not finitely presented for any $n \geq 2$.

\par Let $\Sigma$ be an alphabet and $\Sigma^*$ be the free monoid generated by $\Sigma$, that is, the set of all words with letters in $\Sigma$. For a set of relations $R \subseteq \Sigma^* \times \Sigma^*$ define a \emph{monoid presentation} to be the ordered pair $\lanran{\Sigma \mid R}$, and say that a monoid $\cS$ is \emph{presented} by $\lanran{\Sigma \mid R}$ if $\cS \cong \Sigma^* / \rho_R$ where $\rho_R$ is the smallest congruence on $\Sigma$ containing $R$. 
We say $\cS$ is \emph{finitely presented} if there exists finite $\Sigma$ and finite $R$ such that $\cS$ is presented by $\lanran{\Sigma \mid R}$. 
For $u,v \in \Sigma^*$, we write $u =_\cS v$ to denote that $u$ and $v$ represent the same element of $\cS$.

\subsection[Upper triangular tropical matrices]{A presentation for the monoid of upper triangular tropical matrices}
\par In Section \ref{UTGenSec} we established that $UT_n(\Zmax)$ is finitely generated, we now show that $UT_n(\Zmax)$ is finitely presented for all $n \in \N$, constructing a finite presentation for each $n$ using the minimal generating sets given in Corollary~\ref{UTNZ}.

\par We begin by constructing a finite presentation for $UT_n(\Zmax)$ using a generating set of cardinality $\frac{n(n+5)}{2}$. 
We remark that this is not minimal as the minimal generating set in Corollary~\ref{UTNZ} has cardinality $\frac{n(n+3)}{2} + 1$. 
Nonetheless, this presentation simplifies the normal forms we construct, allowing for more concise proofs. We then use this presentation to construct a presentation with a minimal generating set.

\par First, we define $\Omega_n = \{a_k,a_k^{-1},c_k,d_{ij}\colon 1 \leq k \leq n, \ 1 \leq i < j \leq n \}$, and consider the following relations over $\Omega_n$ for $1 \leq i < j \leq n$ and $1 \leq k, l \leq n$:
\begin{minipage}{.52\linewidth}
\begin{align} 
a_ia_j &= a_ja_i \tag{C1} \label{C1} \\
c_ic_j &= c_jc_i \tag{C2} \label{C2} \\
c_k^2 &= c_k \tag{C3} \label{C3}\\
d_{ij}^2 &= d_{ij} \tag{C4} \label{C4}\\
a_lc_k &= c_ka_l && \tag{C5} \label{C5} \\
a_kd_{ij} &= d_{ij}a_k && i,j \neq k \tag{C6} \label{C6} \\
c_kd_{ij} &= d_{ij}c_k && i,j \neq k \tag{C7} \label{C7}\\
d_{ij}d_{st} &= d_{st}d_{ij} &&  j \neq s < t \neq i \tag{C8} \label{C8} \\
d_{ij}d_{jt} &= d_{jt}d_{ij}d_{it} && j < t \tag{C9} \label{C9} 
\end{align}
\end{minipage}
\begin{minipage}{.41\linewidth}
\begin{align}
d_{ij}a_id_{ij} &= a_id_{ij} \tag{C10} \label{C10}\\
d_{ij}a_jd_{ij} &= d_{ij}a_j \tag{C11} \label{C11}\\
a_ia_jd_{ij} &= d_{ij}a_ia_j \tag{C12} \label{C12}\\
a_kc_k &= c_k \tag{Z1} \label{Z1}\\
c_id_{ij} &= c_i \tag{Z2} \label{Z2}\\
d_{ij}c_j &= c_j \tag{Z3} \label{Z3}\\
a_ka_k^{-1} &= \varepsilon \tag{I1} \label{I1}\\
a_k^{-1}a_k &= \varepsilon \tag{I2}\label{I2} \\
 \nonumber
\end{align}
\end{minipage} 
\newline \\ 
where $\varepsilon$ is the empty word in $\Omega_n^*$. Let $R_n'$ be the collection of all these relations, and note that $a_k$ and $c_k$ commute with all the generators apart from $d_{ik}$ or $d_{kj}$ with $i < k < j$.

We aim to show that $UT_n(\Zmax)$ is presented by $\lanran{\Omega_n \mid R_n'}$ with morphism $\phi \colon \Omega_n^* \to UT_n(\Zmax)$ defined by
\[ \phi(a_k) = A_k(1), \ \phi(a_k^{-1}) = A_k(-1), \ \phi(c_k) = A_k(-\infty), \ \phi(d_{ij}) = E_{ij}\]
for $1 \leq k \leq n$ and $1 \leq i < j \leq n$ and extending multiplicatively.

\par Throughout the rest of this section, we use $\cS$ to denote the monoid presented by $\lanran{\Omega_n \mid R_n'}$. To show that $\cS$ is isomorphic to $UT_n(\Zmax)$, we require two technical lemmas.

\par We begin by showing a number of relations involving $a_k^{-1}$ for $1 \leq k \leq n$ are satisfied by $\cS$. 

\begin{lemma}
The following relations are satisfied by $\cS$. For $1 \leq i < j \leq n$, and $1 \leq k,l \leq n$:

\begin{minipage}{.52\linewidth}
\begin{align} 
a_l^{-1}a_k &= a_ka_l^{-1} \tag{S1} \label{S1}\\
a_i^{-1}a_j^{-1} &= a_j^{-1}a_i^{-1} \tag{S2}\label{S2}\\
a_l^{-1}c_k &= c_ka_l^{-1} && \tag{S3} \label{S3} \\
a_k^{-1}d_{ij} &= d_{ij}a_k^{-1} && i,j \neq k \tag{S4} \label{S4}
\end{align}
\end{minipage}
\begin{minipage}{.41\linewidth}
\begin{align}
a_i^{-1}a_j^{-1}d_{ij} &= d_{ij}a_i^{-1}a_j^{-1} \tag{S5} \label{S5} \\
d_{ij}a_i^{-1}d_{ij} &= d_{ij}a_i^{-1} \tag{S6} \label{S6} \\
d_{ij}a_j^{-1}d_{ij} &= a_j^{-1}d_{ij} \tag{S7} \label{S7} \\
a_k^{-1}c_k &= c_k \tag{S8} \label{S8}
\end{align}
\end{minipage}
\end{lemma}
\begin{proof}
These can be shown by some simple calculations using (\ref{I1}) and (\ref{I2}) with the other relations from $R_n'$. Explicit calculations are given in the author's thesis \cite[Lemma 5.4.1]{TAThesis}.
\end{proof}
% \begin{proof}
% We show that these relations are satisfied by $\cS$ by using (I1) and (I2) with the relations from $R_n'$.
% \begin{enumerate}[(S1):]
%     \item $a_l^{-1}a_k =_S a_l^{-1}a_ka_la_l^{-1} =_S a_l^{-1}a_la_ka_l^{-1} =_S a_ka_l^{-1}$ by (I1), (C1) and (I2).
%     \item $a_i^{-1}a_j^{-1} =_S a_j^{-1}a_ja_i^{-1}a_j^{-1} =_S a_j^{-1}a_i^{-1}a_ja_j^{-1} =_S a_j^{-1}a_i^{-1}$ by (I2), (S1) and (I1).
%     \item $a_l^{-1}c_k =_S a_l^{-1}c_ka_la_l^{-1} =_S a_l^{-1}a_lc_ka_l^{-1} =_S c_ka_l^{-1}$ by (I1), (C5) and (I2).
%     \item $a_k^{-1}d_{i,j} =_S a_k^{-1}d_{i,j}a_ka_k^{-1} =_S a_k^{-1}a_kd_{i,j}a_k^{-1} =_S d_{i,j}a_k^{-1}$ by (I1), (C6) and (I2).
%     \item $a_i^{-1}a_j^{-1}d_{i,j} =_S a_i^{-1}a_j^{-1}d_{i,j}a_ja_ia_i^{-1}a_j^{-1} =_S d_{i,j}a_i^{-1}a_j^{-1}$ by (I1), (C1), (C12) and (I2).
%     \item $d_{i,j}a_i^{-1}d_{i,j} =_S d_{i,j}a_ja_j^{-1}a_i^{-1}d_{i,j} =_S d_{i,j}a_jd_{i,j}a_j^{-1}a_i^{-1} =_S d_{i,j}a_i^{-1}$ by (I1), (S5) and (C11).
%     \item $d_{i,j}a_j^{-1}d_{i,j} =_S d_{i,j}a_j^{-1}a_i^{-1}a_id_{i,j} =_S a_j^{-1}a_i^{-1}d_{i,j}a_id_{i,j} =_S a_j^{-1}d_{i,j}$ by (I2), (S5) and (C10).
%     \item $a_k^{-1}c_k =_S a_k^{-1}a_kc_k =_S c_k$ by (Z1) and (I2).
% \end{enumerate}
% \end{proof}

Note that $a_k^{-1}$ commutes with all the generators except $d_{ik}$ or $d_{kj}$. For each $k \leq n$, let $\Omega_{k,n} = \{a_i,a_i^{-1},c_i,d_{ij} \colon 1 \leq i \leq k, \ 1 \leq i < j \leq n\} \subseteq \Omega_n$ and observe that $\Omega_{n,n} = \Omega_n$. 

This next lemma shows that, given a word over $\Omega_{k,n}$, we can find a word representing the same element with all elements from $\Omega_{k+1,n} \setminus \Omega_{k-1,n}$ to the left of all the elements of $\Omega_{k-1,n}$.

\begin{lemma} \label{aacdpermute}
Let $k \leq n$, $k < h$, and $w \in \Omega_{k-1,n}^*$. Then, $wa_k =_\cS a_kw_1$, $wa_k^{-1}=_\cS a_k^{-1}w_2$, $wc_k =_\cS c_kw_3$, and $wd_{kh} =_\cS d_{kh}w_4$ for some $w_i \in \Omega_{k-1,n}^*$.
\end{lemma}
\begin{proof}
\par Recall that $a_k$, $a_k^{-1}$, and $c_k$ commute with all elements of $\Omega_{k-1,n}$ apart from $d_{i,k}$ with $i < k$. 
Moreover, for all $i < k$,
\begin{align*}
    d_{ik}a_k &=_\cS d_{ik}a_ka_ia_i^{-1} =_\cS a_ka_id_{ik}a_i^{-1} &&\text{by (\ref{I1}), (\ref{C1}), (\ref{C12})}\\
    d_{ik}a_k^{-1} &=_\cS d_{ik}a_k^{-1}a_i^{-1}a_i =_\cS a_k^{-1}a_i^{-1}d_{ik}a_i &&\text{by (\ref{I2}), (\ref{S2}), (\ref{S5})} \\
    d_{ik}c_k &=_\cS c_k &&\text{by (\ref{Z3})}
\end{align*}
Hence, we can permute $a_k$ and $a_k^{-1}$ to the left of $w$, possibly introducing copies of $a_i$ and $a_i^{-1}$ with $i < k$, and we can permute $c_k$ to the left of $w$, removing any $d_{i,k}$ with $i < k$ in $w$. Thus, $wa_k =_\cS a_kw_1$, $wa_k^{-1} =_\cS a_k^{-1}w_2$, and $wc_k =_\cS c_kw_3$ for some $w_1, w_2, w_3 \in \Omega_{k-1,n}$.
\par Finally, note that $d_{kh}$ commutes with all elements of $\Omega_{k-1,n}$ apart from $d_{sk}$ with $s < k$ by (\ref{C6}--\ref{C8}) and (\ref{S4}). However, $d_{sk}d_{kh} =_\cS d_{kh}d_{sk}d_{sh}$ for all $s < k$ by (\ref{C9}).
Hence, we can permute $d_{kh}$ to the left of $w$, possibly introducing some $d_{sh} \in \Omega_{k-1,n}$ with $s < k$. Thus, $wd_{kh} =_\cS d_{kh}w_4$ for some $w_4 \in \Omega_{k-1,n}^*$.
\end{proof}

We are now able to show that $UT_n(\Zmax)$ is finitely presented for all $n \in \N$.
\begin{theorem}
The monoid $UT_n(\Zmax)$ is finitely presented by $\lanran{\Omega_n \mid R_n'}$ for all $n \in \N$.
\end{theorem}
\begin{proof}
We plan to show that $UT_n(\Zmax)$ is isomorphic to $\cS$, the monoid presented by $\lanran{\Omega_n \mid R_n'}$.
For $x \in \Zmax$, $1 \leq k \leq n$, and $1 \leq i < j \leq n$, let 
\[a_k(x) = 
\begin{cases}
a_k^x &x \in \Z, \\
c_i &x = -\infty,
\end{cases} \quad 
d_{i,j}(x) =
\begin{cases}
a_i^xd_{ij}a_i^{-x} &x \in \Z, \\
\varepsilon &x= -\infty,
\end{cases}
\]
and $d_i(x_{i+1},\dots,x_{n}) = d_{i,i+1}(x_{i+1})\cdots d_{i,n}(x_{n})$. 
We aim to show that for any $w \in \Omega_{k,n}^*$, 
\[w =_\cS d_k(x_{k,k+1},\dots x_{k,n})a_k(x_{k,k})v\]
for some $x_{k,j} \in \Zmax$ for $k \leq j \leq n$ and $v \in \Omega_{k-1,n}^*$.

\par Let $w \in \Omega_{k,n}^*$ then, by Lemma~\ref{aacdpermute}, $w =_\cS uv$ for some $u \in (\Omega_{k,n} \setminus \Omega_{k-1,n})^*$ and $v \in \Omega_{k-1,n}^*$. Then, as $\Omega_{k,n} \setminus \Omega_{k-1,n} = \{a_k,a_k^{-1},c_k,d_{kh} \colon k < h \leq n\}$, 
\[ w =_\cS \left( \prod_{i=1}^{\ell'} u_id_{kj_i} \right)u_{\ell'+1}v \]
for some $\ell' \in \N_0$, $k < j_i \leq n$, and $u_i \in \{a_k,a_k^{-1},c_k\}^*$. Since $c_k$ is a zero for $\{a_k,a_k^{-1},c_k\}$ and a left zero for $d_{kh}$ for all $k < h$, it follows that
\[ w =_\cS \left(\prod_{i=1}^{\ell}a_k^{t_i}d_{kj_i}\right)a_k^{t_{\ell+1}}c_k^{\varepsilon_k}v\]
where $\ell \in \N_0$, $t_1,\dots,t_{\ell+1} \in \Z$, $k < j_1,\dots,j_\ell \leq n$, and $\varepsilon_k \in \{0,1\}$.
By the definition of $d_{k,j}(x)$, $a_k^xd_{kj} = d_{k,j}(x)a_k^x$ for $x \in \Z$. Thus,
\[ w =_\cS \left(\prod_{i=1}^{\ell}d_{k,{j_i}}(T_i)\right)a_k^{T_{\ell+1}}c_k^{\varepsilon_k}v.\]
where $T_i = \sum_{j=1}^i t_j$ for $1 \leq i \leq \ell+1$. Now, note that, when $n \neq m$, we can commute the following terms,
\begin{align*}
d_{k,n}(x)d_{k,m}(y) &=_\cS a_k^xd_{kn}a_k^{y-x}d_{km}a_k^{-y} \\
&=_\cS a_k^xd_{kn}a_m^{x-y}a_m^{y-x}a_k^{y-x}d_{km}a_k^{-y} &&\text{(\ref{I1}--\ref{I2})} \\
&=_\cS a_k^xa_m^{x-y}d_{kn}d_{km}a_m^{y-x}a_k^{-x} &&\text{(\ref{C6}), (\ref{C12}), (\ref{S4}--\ref{S5})}\\
&=_\cS a_k^xa_m^{x-y}d_{km}d_{kn}a_m^{y-x}a_k^{-x} &&\text{(\ref{C8})}\\
&=_\cS a_k^yd_{km}a_k^{x-y}a_m^{x-y}a_m^{y-x}d_{kn}a_k^{-x} &&\text{(\ref{C6}), (\ref{C12}), (\ref{S4}--\ref{S5})}\\
&=_\cS a_k^yd_{km}a_k^{x-y}d_{kn}a_k^{-x} &&\text{(\ref{I1}--\ref{I2})}\\
&=_\cS d_{k,m}(y)d_{k,n}(x).
\end{align*} 
When $n = m$, we can simplify in the following way,
\begin{align*}
d_{k,n}(x)d_{k,n}(y) &=_\cS a_k^xd_{kn}a_k^{y-x}d_{kn}a_k^{-y} \\
&=_\cS \begin{cases}
a_k^xd_{kn}(\prod_{i=1}^{|y-x|}a_kd_{kn})a_k^{-y} & y \geq x \\
a_k^xd_{kn}(\prod_{i=1}^{|y-x|}a_k^{-1}d_{kn})a_k^{-y} & y < x
\end{cases} &&\text{(\ref{C10}), (\ref{S6})} \\
&=_\cS \begin{cases}
a_k^xa_k^{y-x}d_{kn}a_k^{-y} & y \geq x \\
a_k^xd_{kn}a_k^{y-x}a_k^{-y} & y < x
\end{cases} &&\text{(\ref{C10}), (\ref{S6})} \\
&=_\cS a_k^{\max(x,y)}d_{kn}a_k^{-\max(x,y)} \\
&=_\cS d_{k,n}(\max(x,y)).
\end{align*} 
\par Now, we define the following variables. For $k < j$, let
\[ 
x_{k,k} = 
\begin{cases}
T_{\ell+1} &\text{if } \varepsilon_k = 0, \\
-\infty &\text{if } \varepsilon_k = 1,
\end{cases} \quad
x_{k,j} = 
\begin{cases}
\max_{j_m = j}(T_m) &\text{if } j_m = j \text{ for some } m, \\
-\infty &\text{otherwise.}
\end{cases}
\]
Then, by above, we have that
\[ w =_\cS d_{k}(x_{k,k+1},\dots,x_{k,n})a_k(x_{k,k})v.\]
Thus, each $w \in \Omega_{k,n}^*$, can be expressed in the above form, and hence by applying this with $k = n,\dots,1$ for $w \in \Omega_n^*$, we get that
\[ w =_\cS a_n(x_{n,n})d_{n-1}(x_{n-1,n}) \cdots a_2(x_{2,2})d_1(x_{1,2},\dots,x_{1,n})a_1(x_{1,1})\]
for some $x_{i,j} \in \Zmax$. 

We can now construct an isomorphism between $\cS$ and $UT_n(\Zmax)$. Define the map $\phi \colon \Omega_n^* \rightarrow UT_n(\Zmax)$, given by $a_i \rightarrow A_i(1), \ a_i^{-1} \mapsto A_i(-1), \ c_i \mapsto A_i(-\infty), \ d_{ij} \mapsto E_{ij}$ and extending multiplicatively. Given $w \in \Omega_n^*$ with the following form 
\[ w = a_n(x_{n,n})d_{n-1}(x_{n-1,n}) \cdots a_2(x_{2,2})d_1(x_{1,2},\dots,x_{1,n})a_1(x_{1,1}),\] 
a simple calculation shows that
\[ \phi(w) =
\begin{pmatrix}
x_{1,1} & \dots & x_{1,n} \\
& \ddots & \vdots \\
& & x_{n,n}
\end{pmatrix}.\]
Thus, as $x_{i,j} \in \Zmax$ is arbitrary for all $i \leq j$, every matrix in $UT_n(\Zmax)$ is the image of a word of the above form, and hence the set of words of the above form are in bijection with $UT_n(\Zmax)$. So, it suffices to check that the images of the generators satisfy the relations in $R_n'$. These are simple calculations that can be easily checked, but for the explicit calculations, see the author's thesis \cite[Theorem 5.4.4]{TAThesis}.
Thus, $UT_n(\Zmax)$ is finitely presented by $\lanran{\Omega_n \ | \ R_n'}$.
\end{proof}

If a semigroup is finitely presented then it can be finitely presented with every finite generating set for the semigroup \cite[Proposition 3.1]{RuskucThesis}.
So, we use the above theorem to construct a finite presentation for $UT_n(\Zmax)$ using the minimal and irreducible generating set from Corollary~\ref{UTNZ}. For this, we define the alphabet $\Sigma_n = \{a_k,b,c_k,d_{ij}\colon 1 \leq k \leq n, \ 1 \leq i < j \leq n \}$ and the relations, for $1 \leq i < j \leq n$, and $1 \leq k \leq n$:

\begin{minipage}{0.33\linewidth}
\begin{equation}
    a_kb = ba_k, \tag{R1} \label{R1} 
\end{equation}
\end{minipage}
\begin{minipage}{0.33\linewidth}
\begin{equation}
    d_{ij}b = bd_{ij}, \tag{R2} \label{R2}
\end{equation}
\end{minipage}
\begin{minipage}{0.33\linewidth}
\begin{equation}
    a_1\cdots a_nb = \varepsilon \tag{R3} \label{R3}.
\end{equation}
\end{minipage} 
\newline 
\par We define $R_n$ to be the collection of relations (\ref{C1}--\ref{C11}), (\ref{Z1}--\ref{Z3}), and (\ref{R1}--\ref{R3}). That is $R_n'$ with (\ref{I1}--\ref{I2}), and (\ref{C12}) replaced with (\ref{R1}--\ref{R3}).

\begin{theorem}
The monoid $UT_n(\Zmax)$ is finitely presented by $\lanran{\Sigma_n \mid R_n}$ for all $n \in \N$.
\end{theorem}
\begin{proof}
Let $\cM$ be the monoid presented by $\lanran{\Sigma_n \mid R_n}$ and recall that $\cS \cong UT_n(\Zmax)$ is the monoid presented by $\lanran{\Omega_n \mid R_n'}$. We show that $\cM \cong \cS$. Define $\phi \colon \cM \rightarrow \cS$ to be the map given by $a_i \mapsto a_i, \ c_i \mapsto c_i, \ d_{ij} \mapsto d_{ij},$ and $b \mapsto a_1^{-1} \cdots a_n^{-1}$ and extending multiplicatively. 
To see that $\phi$ is a well-defined map, we must show that $\phi(\Sigma_n^*)$ satisfies the relations $R_n$. So, note that $\phi$ is the identity map on $\Sigma_n \setminus \{b\}$, and hence satisfies the relations (\ref{C1}--\ref{C11}) and (\ref{Z1}--\ref{Z3}). Thus, it suffices to check that $\phi(\Sigma_n^*)$ satisfies the relations (\ref{R1}--\ref{R3}). For $1 \leq i < j \leq n$ and $1 \leq k \leq n$,
\begin{align*}
\phi(a_k)\phi(b) &= a_ka_1^{-1} \cdots a_n^{-1} =_\cS a_1^{-1} \cdots a_n^{-1}a_k = \phi(b)\phi(a_k) &&\text{by (\ref{S1}),} \\
\phi(d_{ij})\phi(b) &= d_{ij}a_1^{-1} \cdots a_n^{-1} =_\cS a_1^{-1} \cdots a_n^{-1}d_{ij} = \phi(b)\phi(d_{ij}) &&\text{by (\ref{S2}), (\ref{S4}--\ref{S5}),} \\
\phi(a_1)\cdots &\phi(a_n)\phi(b) = a_1 \cdots a_na_1^{-1}\cdots a_n^{-1} =_\cS \varepsilon = \phi(\varepsilon) &&\text{by (\ref{C1}), (\ref{I1}).}
\end{align*}
\par Now, define $\psi \colon \cS \rightarrow \cM$ to be the map given by $a_i \mapsto a_i, \ c_i \mapsto c_i, \ d_{ij} \mapsto d_{ij},$ and $a_i^{-1} \mapsto a_1 \cdots a_{i-1}a_{i+1}\cdots a_nb$ and extending multiplicatively. 
To show that $\psi$ is a well-defined map, we show that $\psi(\Omega_n^*)$ satisfies the relations $R_n'$. 
Again, note that $\psi$ is the identity map on $\Omega_n \setminus \{a_k^{-1} \colon 1 \leq k \leq n\}$, and hence satisfies the relations (\ref{C1}--\ref{C11}) and (\ref{Z1}--\ref{Z3}). Thus, it suffices to check that $\psi(\Omega_n^*)$ satisfies the relations (\ref{I1}), (\ref{I2}), and (\ref{C12}). For $1 \leq i < j \leq n$ and $1 \leq k \leq n$,

\begin{align*}
\psi(a_k)\psi(a_k^{-1}) &= a_ka_1 \cdots a_{k-1}a_{k+1}\cdots a_nb =_\cM \varepsilon &&\text{(\ref{C1}), (\ref{R3}),} \\
\psi(a_k^{-1})\psi(a_k) &= a_1 \cdots a_{k-1}a_{k+1}\cdots a_nba_k =_\cM \varepsilon &&\text{(\ref{C1}), (\ref{R1}), (\ref{R3}),} \\
\psi(a_i)\psi(a_j)\psi(d_{ij}) &= a_ia_jd_{ij} \\ 
&=_\cM a_ia_jd_{ij}a_1 \cdots a_nb &&\text{(\ref{R3}),}\\
&=_\cM a_1 \cdots a_nbd_{ij}a_ia_j &&\text{(\ref{C1}), (\ref{C6}), (\ref{R1}--\ref{R2}),}\\
&=_\cM d_{ij}a_ia_j &&\text{(\ref{R3}),} \\ 
&= \psi(d_{ij})\psi(a_i)\psi(a_j).
\end{align*}
Thus, $\phi$ and $\psi$ are well defined morphisms. To see that $\phi$ and $\psi$ are mutually inverse morphisms, note that $\psi \phi(a_i) = a_i$, $\psi \phi(c_i) = c_i$, $\psi \phi(d_{ij}) = d_{ij}$, and
\begin{align*}
\psi \phi(b) &= \psi(a_1^{-1} \cdots a_n^{-1}) \\
&= (ba_2\cdots a_n)\cdots (ba_1 \cdots a_{n-1}) \\
&=_\cM a_1^{n-1}\cdots a_n^{n-1}b^n &&\text{by (\ref{C1}), (\ref{R1}),} \\
&=_\cM (a_1\cdots a_nb)^{n-1}b &&\text{by (\ref{C1}), (\ref{R1}),} \\
&=_\cM b &&\text{by (\ref{R3}).} 
\end{align*}
Therefore, $\psi\phi \colon \cM \rightarrow \cM$ is the identity map on $\cM$. Similarly, $\phi \psi(a_i) = a_i, \phi \psi(c_i) = c_i$, and $\phi \psi(d_{ij}) = d_{ij}$, so finally note that, for $1 \leq k \leq n$,
\begin{align*}
\phi \psi(a_k^{-1}) &= \phi(a_1\cdots a_{k-1}a_{k+1}\cdots a_nb) \\
&= a_1\cdots a_{k-1}a_{k+1}\cdots a_na_1^{-1}\cdots a_n^{-1} \\
&=a_k^{-1} &&\text{by (\ref{C1}), (\ref{S1}), (\ref{I1}).}
\end{align*}
Thus, $\phi\psi \colon \cS \rightarrow \cS$ is the identity map on $\cS$. Therefore, $\phi$ and $\psi$ are mutually inverse morphisms and $\cM$ and $\cS$ are isomorphic. Hence, $\lanran{\Sigma_n \mid R_n}$ is a finite presentation for $UT_n(\Zmax)$ with a minimal generating set.
\end{proof}

\begin{remark}
    The presentation for $UT_n(\Zmax)$ given in the above theorem has $\frac{n(n+3)}{2}+1$ generators and $\frac{1}{8}(n^4+6n^3+15n^2+10n+8)$ relations.
\end{remark}

\subsection{Full tropical matrix monoids presentations}
In Section~\ref{FullGenSec}, we showed that $M_n(\Zmax)$ is finitely generated if and only if $n \leq 2$. We now show that $M_2(\Zmax)$ is not finitely presented. But, first, we need the following result which describes the $\cJ$-class structure of $M_2(\Zmax)$.

\begin{proposition}[{\cite[Corollary 3.8]{JK2by2Structure}}] \label{prop:M2JClasses}
    The $\cJ$-classes of $M_2(\Zmax)$ are linearly ordered indexed by the set $\cI = \{-\infty\} \cup \N_0 \cup \{\infty, \infty^*\}$ where $-\infty < n < \infty < \infty^*$ for all $n \in \N_0$. Moreover, the $\cJ$-classes can be expressed in the following way, for $n \in \N_0$,
    \begin{align*}
        J_{\infty^*} &= GL_2(\Zmax), \\
        J_{\infty} &= \{A \colon A \text{ has exactly one } -\infty \text{ entry} \}, \\
        J_n &= \{A \colon A_{11}+A_{22} = A_{12} +A_{21} + n \text{ or } A_{11}+A_{22} +n = A_{12} +A_{21} \}, \\
        J_{-\infty} &= \left\{\begin{pmatrix}
            -\infty &-\infty \\
            -\infty &-\infty
        \end{pmatrix}\right\}.
    \end{align*}
\end{proposition}
\begin{theorem} \label{thm:2by2NotFinPres}
$M_2(\Zmax)$ is not finitely presented.
\end{theorem}
\begin{proof}
\par Let $X = \{a,a^{-1},b,c,d\}$ and  $\sigma \colon X^* \rightarrow M_2(\Zmax)$ be the morphism obtained by mapping 
$\sigma(a) = A_1(1)$, $\sigma(a^{-1}) = A_1(-1)$,  $\sigma(b) = \left( \begin{smallmatrix}
    -\infty & 0 \\
    0 & -\infty
\end{smallmatrix}\right)$, $\sigma(c) = A_1(-\infty)$, $\sigma(d) = E_{12}$ 
and extending multiplicatively. 
Note that $\sigma$ is surjective as these matrices form a generating set for $M_2(\Zmax)$ as they generate the matrices in Corollary~\ref{2by2Z}. Now, let $\lanran{X \mid R}$ be a presentation for $M_2(\Zmax)$ and, for a contradiction, suppose that $R$ is finite. 

By Proposition~\ref{prop:M2JClasses}, the $\cJ$-classes of $M_2(\Zmax)$ are linearly ordered with index set $\cI = \{ -\infty\} \cup \N_0 \cup \{\infty, \infty^*\}$ where $-\infty < n < \infty < \infty^*$ for all $n \in \N_0$. 

\par Let $R = R_1 \cup R_2$ where $R_1$ is the set of relations $(u,v) \in R$ with $\sigma(u) \in J_{\infty^*} \cup J_{\infty}$ and $R_2$ is the set of relations $(u,v) \in R$ with $\sigma(u) \in J_r$ for some $r \in \N_0 \cup \{-\infty\}$. Moreover, as there are only finitely many relations there exists $n \in \N_0$ such that $\sigma(u) \notin \bigcup_{i \geq  n} J_i$ for all $(u,v) \in R_2$. 

\par Next, consider the word $da^nbd$, a simple calculation gives that
\[ \sigma(da^nbd) = 
\begin{pmatrix}
    0 & 0 \\
      & 0
\end{pmatrix}
\begin{pmatrix}
    1 &  \\
      & 0
\end{pmatrix}^n
\begin{pmatrix}
      & 0 \\
    0 & 
\end{pmatrix}
\begin{pmatrix}
    0 & 0 \\
      & 0
\end{pmatrix}
= 
\begin{pmatrix}
    0 & n \\
    0 & 0
\end{pmatrix}.
\]
Remark that $\sigma(da^nbd) \in J_n$ by Proposition~\ref{prop:M2JClasses}. Now, let $w = a^nbda^nbd$ and $w' = da^nbda^nb$ and note that
\[\sigma(w) = 
\begin{pmatrix}
     & n \\
    0 & 
\end{pmatrix}
\begin{pmatrix}
    0 & n \\
    0 & 0
\end{pmatrix}
=
\begin{pmatrix}
    n & n \\
    0 & n
\end{pmatrix}
=
\begin{pmatrix}
    0 & n \\
    0 & 0
\end{pmatrix}
\begin{pmatrix}
      & n \\
    0 & 
\end{pmatrix}
= \sigma(w').\]
By Proposition~\ref{prop:M2JClasses}, $\sigma(w), \sigma(w') \in J_n$, and hence the equality $\sigma(w) = \sigma(w')$ is a consequence only of the relations in $R$ corresponding to elements in the $\cJ$-classes in and above $J_n$, that is, exactly the relations in $R_1$. 

\par Note that $w$ contains exactly one $b$ to the left of the first $d$ while $w'$ does not and $\sigma(da^nd) \in J_\infty$ while $\sigma(da^nbd) \in J_n$. 
Thus, to turn $w'$ into $w$, it is necessary that $\sigma(da^nb) = \sigma(s)$ for some $s \in X^*$ with at least one $d$ and exactly one $b$ before the first $d$. 
Moreover, if $p \in \{a,a^{-1},b\}^*$ contains $b$ an odd number of times, then  $\sigma(dpd) \in J_i$ for some $i \in \N_0$. Thus, $b$ occurs an even number of times between any two occurrences of $d$ in $s$.

\par Recall $\phi_2 \colon M_2(\Zmax) \to M_2(\B)$ is the morphism mapping the entries from $\Z$ to $1$ and the $-\infty$ entries to $0$. 
Thus, $\phi_2(\sigma(da^nb)) = \phi_2(\sigma(db))$ and $\phi_2(\sigma(s)) = \phi_2(\sigma(bdb^\varepsilon))$ for some $\varepsilon \in \{0,1\}$ as $\phi_2(A_1(1)) = \phi_2(A_1(-1)) = B^2 = I_2$, and $E_{12}^2 = E_{12}$. 
Then, as $\sigma(da^nb) = \sigma(s)$, it follows that $\phi_2(\sigma(da^nb)) = \phi_2(\sigma(s))$. 

However, $\phi_2(\sigma(db)) \neq \phi_2(\sigma(bdb^\varepsilon))$ for either $\varepsilon \in \{0,1\}$. Thus, the relation $(a^nbda^nbd,da^nbda^nb)$ is not implied by the relations in $R$ giving a contradiction, so $R$ is infinite. Therefore, $M_2(\Zmax)$ is not finitely presented, as a finitely presented semigroup can be finitely presented with any finite generating set for the semigroup \cite[Proposition 3.1]{RuskucThesis}.
\end{proof}

\begin{corollary}
    The monoid $M_n(\Zmax)$ is finitely presented if and only if $n = 1$.
\end{corollary}
\begin{proof}
    It is clear $M_1(\Zmax) \cong (\Zmax, \cdot)$ is finitely presented, while $M_n(\Zmax)$ is not finitely presented for $n \geq 2$ by Corollary~\ref{2by2Z} and Theorem~\ref{thm:2by2NotFinPres}.
\end{proof}

\section{Growth of commutative bipotent matrices} \label{GrowthSec}

In this section, we establish upper bounds for the growth of finitely generated subsemigroup of $M_n(S)$ or $UT_n(S)$, when $S$ is a commutative bipotent semiring. When $S = \T$, we produce more explicit bounds, depending only on $n$ and the rank of the free abelian group generated by the finite matrix entries. Finally, we show that the bounds of the polynomial degree are sharp by producing a family of examples.

\par For a semigroup $S$ generated by the finite set $X$, the \emph{growth function of $S$ with respect to} $X$ is $f_X(k) = |\cup^{k}_{i=1}X^i|$. We say $f_X(k)$ is \emph{bounded above (resp. below) by a polynomial of degree} $n$ if there exists $c_X > 0$ such that for all $k \in \N$, $f_X(k) \leq c_Xk^n$ (resp. $f_X(k) \geq c_Xk^n$).
\par It is well-known that if the growth function of $S$ with respect to $X$ is bounded above (resp. below) by a polynomial of degree $n$ then the growth function with respect to any finite generating set is bounded above (resp. below) by a polynomial of degree $n$. So, we may say that $S$ has growth function bounded above/below by a polynomial of degree $n$, without reference to a generating set.

\subsection{Upper bounds for growth} \label{GrowthSecSec}

We say that a semiring $(S,+,\cdot)$ is \emph{bipotent} if $x + y \in \{x,y\}$ for all $x,y \in S$. We begin by finding upper bounds for the growth of finitely generated subsemigroup of $M_n(S)$ or $UT_n(S)$ when $S$ is a bipotent semiring.

\begin{proposition} \label{prop:fullmatrixgrowth}
Let $S$ be a bipotent semiring, $X \subseteq M_n(S)$ be a finite set, and $T = \lanran{X}$. If the growth of the multiplicative semigroup generated by the entries of the matrices in $X$ is bounded above by a polynomial of degree $t \in \N_0$. Then, the growth function of $T$ is bounded above by a polynomial of degree $tn^2$.
\end{proposition}
\begin{proof}
For every $k \geq 1$, let $C_k$ be the set of all the non $0_S$ entries of the matrices in $X^k$ and let $c_k = |\cup_{i=1}^k C_i|$.
As the growth of the semigroup generated by the entries of the matrices in $X$ is bounded above by a polynomial of degree $t$, we have that $c_k \leq \beta k^t$ for some $\beta > 0$ as $S$ is bipotent.
\par Hence, as every matrix in $X^k$ has entries in $C_k \cup \{0_S\}$, we obtain for every $k \in \N$,
\[ f_X(k) \leq (c_k+1)^{n^2} \leq  (\beta k^{t}+1)^{n^2} \leq ((\beta+1) k^{t})^{n^2} = \delta k^{tn^2} \]
where $\delta = (\beta+1)^{n^2}$.
\end{proof}
\par If $S$ is a commutative bipotent semiring and $X \subseteq S$ with $|X| = t$. Then, the subsemigroup of the multiplicative semigroup $(S,\cdot)$ generated by $X$ has growth bounded above by a polynomial of degree $t$ as it is a quotient of the free commutative semigroup $\N^t$ which has growth bounded above by a polynomial of degree $t$.
Thus, we may apply the above theorem to any finitely generated subsemigroup of $M_n(S)$ when $S$ is a commutative bipotent semiring.

\par The upper bound on the degree of the growth for $M_n(S)$ where $S$ is a commutative bipotent semiring given in \cite{CombinatorialSemigroups} is $(c-1)n^2 + 1$, where $c$ is the number of distinct matrix entries in the generating set $X$. 
Thus, the new bound given above is only worse when $n \geq 2$ and the growth of the multiplicative semigroup generated by the entries of the matrices in $X$ is bounded below by a polynomial of degree $c$. In particular, no matrix in $X$ has $0_S$ or $1_S$ as an entry.

\par To achieve a more explicit upper bound of the polynomial degree, we restrict to the case where $S = \T$. But, we first require the following lemma which gives the growth of finitely generated subsemigroups of the multiplicative semigroup of $\T$ in terms of the rank of free abelian subgroup they generate as a group.

\begin{lemma} \label{rbound}
Let $C \subseteq (\R \cup \{-\infty\},+)$ be a finite set and $T = \lanran{C}$. Then the growth of $T$ is bounded above by a polynomial of degree $t$, where $t$ is the rank of the free abelian group generated, as a group, by $C \setminus \{-\infty\}$. 
\end{lemma}
\begin{proof}
Let $G$ be the free abelian group generated, as a group, by $D = C \setminus \{-\infty\}$. Let $f(k)$ be the growth of $T$ with respect to $C$ and $g(k)$ be the growth of $G$ with respect to $(D \cup D^{-1})$, where $D^{-1} = \{ d^{-1} \colon d \in D\}$. 
Clearly, $f(k) \leq g(k) + 1$ as $-\infty \notin D$. 
Moreover, as $G$ is a free abelian group of rank $t$, $G$ has growth upper bounded by a polynomial of degree $t$ \cite[Theorem 3.2]{WolfZGrowth}.
Thus, 
\[f(k) \leq g(k) + 1 \leq ck^t\] 
for some $c > 0$.
\end{proof}

\begin{corollary} \label{cor:fullmatrixT}
Let $T$ be a finitely generated subsemigroup of $M_n(\T)$ and $t$ be the rank of the free abelian subgroup of $(\R,+)$ generated as a group by the finite entries of the matrices in $T$. Then, the growth function of $T$ is bounded above by a polynomial of degree $tn^2$.
\end{corollary}
\begin{proof}
The finite entries of the matrices in $T$ and the finite entries of the matrices in any generating set for $T$ generate, as a group, the same free abelian subgroup of $(\R,+)$. Thus, as $\T$ is bipotent, the result follows immediately from Proposition~\ref{prop:fullmatrixgrowth} and Lemma~\ref{rbound} 
\end{proof}

\par If we consider the case where $S = \Qmax$, then we can further simplify the result.

\begin{corollary}
Let $T$ be a finitely generated subsemigroup of $M_n(\Qmax)$. Then, the growth function of $T$ is polynomially upper bounded of degree $n^2$.
\end{corollary}
\begin{proof}
All finitely generated subgroups of $(\Q,+)$ are either trivial or isomorphic to $(\Z,+)$, \cite[Exercise 4.2.6]{SubgroupsQ}.
\end{proof}
\par We now provide similar results for the semigroup of upper triangular matrices over bipotent semirings.

\begin{proposition} \label{uppertri}
Let $S$ be a bipotent semiring, $X \subseteq M_n(S)$ be a finite set, and $T = \lanran{X}$. If the growth of the multiplicative semigroup generated by the entries of the matrices in $X$ is bounded above by a polynomial of degree $t \in \N_0$. Then, the growth function of $T$ is bounded above by a polynomial of degree $\frac{tn(n+1)}{2}$.
\end{proposition}
\begin{proof}
    Identical to the proof of Proposition~\ref{prop:fullmatrixgrowth}
\end{proof}

\par Now, we again restrict to the cases where the bipotent semiring is $\T$ or $\Qmax$ to give explicit bounds on the growth of finitely generated subsemigroups of $UT_n(\T)$ and $UT_n(\Qmax)$.

\begin{corollary}\label{uppertriT}
Let $T$ be a finitely generated subsemigroup of $UT_n(\T)$ and $t$ be the rank of the free abelian subgroup of $(\R,+)$ generated as a group by the finite entries of the matrices in $T$. Then, the growth function of $T$ is bounded above by a polynomial of degree $\frac{tn(n+1)}{2}$.
\end{corollary} 
\begin{proof}
    Identical to the proof of Corollary~\ref{cor:fullmatrixT} using Lemma~\ref{rbound} and Proposition~\ref{uppertri}.
\end{proof}

\begin{corollary}
Let $T$ be a finitely generated subsemigroup of $UT_n(\Qmax)$. Then, the growth function of $T$ is bounded above by a polynomial of degree $\frac{n(n+1)}{2}$.
\end{corollary}
\begin{proof}
All finitely generated subgroups of $(\Q,+)$ are either trivial or isomorphic to $(\Z,+)$, \cite[Exercise 4.2.6]{SubgroupsQ}.
\end{proof}

\subsection{The bounds are sharp}
We now show that for all $n \in \N$ and $t \in \N_0$, there exist finitely generated subsemigroups of $M_n(\T)$ and $UT_n(\T)$ such that the finite entries generate, as a group, a free abelian group of rank $t$ and, the growth functions are bounded below by polynomials of degrees $tn^2$ and $\frac{tn(n+1)}{2}$ respectively, that is, the upper bounds given by Corollary~\ref{cor:fullmatrixT} and Corollary~\ref{uppertriT}.

\begin{theorem} \label{thm:UpperTriGrowthSharp}
Let $n \in \N$ and $t \in \N_0$. Then, there exists a finite set $X \subseteq UT_n(\T)$ such that the growth function of $\lanran{X}$ is bounded below by $ck^{\frac{tn(n+1)}{2}}$ for some $c > 0$ where $t$ is the rank of the free abelian subgroup of $(\R,+)$ generated as a group by the finite entries of the matrices in $\lanran{X}$.
\end{theorem}
\begin{proof}
The proof is immediate if $t = 0$, so we may assume $t \geq 1$.
Let $I = \{\gamma_1,\dots, \gamma_t\} \subseteq (\R,+)$ be a minimal group generating set for a free abelian group of rank $t$. 
Consider the set of matrices $\cM_k \subseteq UT_n(\T)$ such that the entries on and above the diagonal are the tropical product of at most $\floor{\frac{k-n}{2n}}$ elements from $I$. 
Now, let $X$ be the set of all $n \times n$ upper triangular matrices with entries from $\{ \gamma_1,\dots, \gamma_t, -\gamma_1, \dots, -\gamma_t,0,-\infty\}$. We now show that $\cM_k \subseteq X^k$ for all $k \in \N$.

\par Let $A \in \cM_k$ and $L_m,R_m \in UT_n(\T)$ be diagonal with $(L_m)_{ii} = A_{im}$ for all $i \leq m$ and 0 otherwise and $(R_m)_{ii} = -A_{im}$ for all $i < m$ and 0 otherwise. Let $E'_m \in UT_n(\T)$ with all diagonal entries being 0, $(E'_m)_{im} = 0$ for all $i \leq m$ and $-\infty$ otherwise. To show that $A \in X^k$, let
\[ \Sigma = \prod^{n-1}_{m=0} L_{n-m}E'_{n-m}R_{n-m} \]
and note that $(L_mE'_mR_m)_{ii} = 0$ for $i \neq m$. Thus, for $i \leq j$,
\[ \Sigma_{ij} = (L_j)_{ii}(E'_j)_{ij}(R_j)_{jj} = A_{ij} + 0 + 0 = A_{ij},\]
and hence, $A = \Sigma$.
Then, as $A \in \cM_k$, each $A_{ij}$ can be expressed as the product of at most $\floor{\frac{k-n}{2n}}$ entries from $I$, so the diagonal matrices $L_m$ and $R_m$ can be expressed as the product of at most $\floor{\frac{k-n}{2n}}$ matrices from $X$.
Therefore, for each $1 \leq m \leq n$, $L_mE'_mR_m$ can be expressed as the product of $2\floor{\frac{k-n}{2n}}+1 \leq \floor{\frac{k}{n}}$ matrices from $X$.
Hence, $A$ can be expressed as the product of $n\floor{\frac{k}{n}} \leq k$ matrices from $X$, and thus $\cM_k \subseteq X^k$.

\par Now, as $\{\gamma_1,\dots,\gamma_t\}$ is a minimal group generating set for a free abelian group, the monoid generated by $\{\gamma_1,\dots,\gamma_t\}$ is a free commutative monoid of rank $t$ and has a growth function bounded below by $c'k^t$ for some $c' > 0$. Thus, there are at least $(c'(\floor{\frac{k-n}{2n}})^t)^{\frac{n(n+1)}{2}}$ matrices in $\cM_k$. 

\par Therefore, there exists $c > 0$ such that $|\cM_k| + |X|$, and hence the growth function of $X$, is bounded below by $ck^{\frac{tn(n+1)}{2}}$.
\end{proof}

\begin{theorem} \label{thm:FullGrowthSharp}
Let $n \in \N$ and $t \in \N_0$. Then, there exists a finite set $X \subseteq M_n(\T)$ such that the growth function of $\lanran{X}$ is bounded below by $ck^{tn^2}$ for some $c > 0$ where $t$ is the rank of the free abelian subgroup of $(\R,+)$ generated as a group by the finite entries of the matrices in $\lanran{X}$.
\end{theorem}
\begin{proof}
The proof is immediate if $t = 0$, so we may assume $t \geq 1$.
Let $I = \{\gamma_1,\dots, \gamma_t\} \subseteq (\R,+)$ be a minimal group generating set for a free abelian group of rank $t$ such that $1 \leq \gamma_i \leq 2$ for each $i$.
Now, consider the set of matrices $\cM_k \subseteq M_n(\T)$ such that the diagonal entries of the matrices are the tropical product of between $\floor{\frac{2k-4n}{16n+3}}$ and $\floor{\frac{3k-6n}{16n+3}}$ elements from $I$ and the off-diagonal entries are the tropical product of between $0$ and $\floor{\frac{k-2n}{16n+3}}$ elements from $-I = \{-\gamma_1,\dots,-\gamma_t\}$.
Let $X$ be the set of all $n \times n$ matrices with entries from $I \cup -I \cup \{0, -\infty\}$.
We now show that $\cM_k \subseteq X^k$ for all $k \in \N$.
\par Let $A \in \cM_k$ and $L_m,R_m \in M_n(\T)$ be the diagonal matrix with entries $(L_m)_{ii} = A_{im} - A_{mm}$ if $i > m$ and 0 otherwise and $(R_m)_{ii} = A_{im} - A_{ii}$ if $i < m$ and 0 otherwise. 
Let $E_m$ be the matrix where all diagonal entries are 0, $(E_m)_{im} = 0$ for all $i \geq m$, and all other entries are $-\infty$. 
Similarly, let $E'_m$ be the matrix where all diagonal entries are 0, $(E'_m)_{im} = 0$ for all $i \leq m$, and all other entries are $-\infty$. Let $\Lambda$ be the diagonal matrix where $\Lambda_{ii} = A_{ii}$ for all $1 \leq i \leq n$.
To show that $A \in X^k$, let
\[
\Sigma = \left(\prod_{m=1}^{n} L_mE_mL_{m}^{-1}\right)\Lambda \left(\prod_{m=0}^{n-1} R_{n-m}E'_{n-m}R_{n-m}^{-1}\right)
\]
and note that $(L_mE_mL_{m}^{-1})_{ii} =(R_mE'_mR_{m}^{-1})_{ii} = 0$.
Thus, for $1 \leq i,j \leq n$,
\[ \Sigma_{ij} = \max_{m \leq i,j} ((L_m)_{ii} + \Lambda_{mm} + (R_j)_{mm}) = \max_{m \leq i,j} (A_{im} + A_{mj} - A_{mm}) = A_{ij}\]
as if $m < \min(i,j)$ then $A_{im} + A_{mj} - A_{mm} \leq 0 + 0 - \floor{\frac{2k-4n}{16n+3}} \leq A_{ij}$ as $1 \leq \gamma_s \leq 2$ for each $s$. Thus, we have that $\Sigma = A$. 
\par Now, for all $1 \leq m \leq n$, $E_m, E_m' \in X$ and both $L_m$ and $R_m$ (and therefore also $L_m^{-1}$ and $R_m^{-1})$ can be expressed as the product of $\floor{\frac{3k-6n}{16n+3}}+ \floor{\frac{k-2n}{16n+3}}$ matrices from $X$. Similarly, $\Lambda$ can be expressed as the product of $\floor{\frac{3k-6n}{16n+3}}$ matrices from $X$. Thus, $\Sigma$ can be expressed as the product of
\begin{multline*}
    4n\floor{\frac{3k-6n}{16n+3}}+4n\floor{\frac{k-2n}{16n+3}}+2n + \floor{\frac{3k-6n}{16n+3}} \\
    \leq \frac{16n(k-2n)}{16n+3}+2n + \frac{3k-6n}{16n+3} 
= k
\end{multline*}
matrices from $X$, and hence $A \in X^k$. 
\par Now, as $I$ is a minimal group generating set for a free abelian group, the monoid generated by $I$ is a free commutative monoid of rank $t$ and has a growth function bounded below by $c'k^t$ for some $c' > 0$. Thus, there are at least $(c'(\floor{\frac{k-2n}{16n+3}})^t)^{n^2}$ matrices in $\cM_k$. 

Therefore, there exists a $c > 0$ such that $|\cM_k| + |X|$, and hence the growth function of $X$, is bounded below by $ck^{tn^2}$.
\end{proof}

From the Theorems~\ref{thm:UpperTriGrowthSharp} and \ref{thm:FullGrowthSharp}, we have now shown that the bounds given in Corollary~\ref{cor:fullmatrixT} and \ref{uppertriT} are sharp.
\begin{corollary}
For all $n \in \N$ and $t \in \N_0$, there exist finitely generated subsemigroups of $M_n(\T)$ and $UT_n(\T)$ such that their growth functions are bounded above and below by polynomials of degree $tn^2$ and $\frac{tn(n+1)}{2}$ respectively where $t$ is the rank of the free abelian group generated as a group by the finite entries of the matrices in subsemigroup.
\end{corollary}

\par \textbf{Acknowledgements.} The author thanks Marianne Johnson and Mark Kambites for many helpful conversations.

\bibliographystyle{plain}
\bibliography{Bib.bib}
\end{document}